\documentclass[a4paper]{amsart}
\usepackage[hidelinks]{hyperref}
\usepackage[usenames,dvipsnames]{color}
\usepackage{amssymb, amsthm}
\usepackage{mathrsfs}
\usepackage{cleveref}
\usepackage{tikz-cd}
\usepackage{bm}
\usepackage{enumerate}
\usepackage{ulem}

\newtheorem{theorem}[equation]{Theorem}
\newtheorem*{theorem*}{Theorem}
\newtheorem{proposition}[equation]{Proposition}
\newtheorem{corollary}[equation]{Corollary}
\newtheorem{lemma}[equation]{Lemma}
\newtheorem{construction}[equation]{Construction}
\theoremstyle{definition}
\newtheorem{definition}[equation]{Definition}
\theoremstyle{remark}
\newtheorem{remark}[equation]{Remark}
\newtheorem{example}[equation]{Example}

\DeclareMathOperator{\finitetype}{ftype}

\newcommand{\finp}{\fin_{\OK}^{p}}
\newcommand{\finpN}{\fin_{\OK}^{p,N}}
\newcommand{\ExtMon}{\Ext_{\fin^{p}_{\OK}}^{\nu}}
\newcommand{\ExtMonR}{\Ext_{\fin^{p}_{R}}^{\nu}}
\newcommand{\ShAbfppf}[1]{\Sh_{\Ab}(\Sch_{#1, \fppf})}
\newcommand{\sQo}{\sQ^{1}_{\varpi}}
\newcommand{\sQon}{\sQo[n]}
\newcommand{\Resoo}{\Res^{\Gamma_{K}}_{\Gamma_{\varpi^{\flat}}}}
\newcommand{\Reptors}{\Rep_{\tors}}
\newcommand{\ModZpGK}{\Mod_{\Zp,\Gamma_{K}}}
\newcommand{\ModtorsZpGK}{\Mod_{\Zp,\Gamma_{K}}^{\tors}}
\newcommand{\RepQp}{\Rep_{\Qp}(\Gamma_{K})}
\newcommand{\RepZp}{\Rep_{\Zp}(\Gamma_{K})}
\newcommand{\RepcrysrQp}{\Rep^{\crys,r}_{\Qp}(\Gamma_{K})}
\newcommand{\Repcrysrfree}{\Rep^{\crys,r}_{\free}(\Gamma_{K})}
\newcommand{\Repcrysrtors}{\Rep^{\crys,r}_{\tors}(\Gamma_{K})}
\newcommand{\Repcrysrft}{\Rep^{\crys, r}_{\finitetype}(\Gamma_{K})}
\newcommand{\ModcrysrZpGK}{\Mod^{\crys,r}_{\Zp,\Gamma_{K}}}
\newcommand{\RepcrysoneZp}{\Rep^{\crys,1}_{\free}(\Gamma_{K})}
\newcommand{\Repcrysonetors}{\Rep^{\crys,1}_{\tors}(\Gamma_{K})}
\newcommand{\Reponesst}{\Rep^{1, \stbl}_{\Zp}(\Gamma_{K})}
\newcommand{\Reponessttors}{\Rep^{1, \stbl}_{\tors}(\Gamma_{K})}
\newcommand{\Repsstrfree}{\Rep^{\stbl,r}_{\free}(\Gamma_{K})}
\newcommand{\Repsstrtors}{\Rep^{\stbl,r}_{\tors}(\Gamma_{K})}
\newcommand{\Repsstonetors}{\Rep^{\stbl,1}_{\tors}(\Gamma_{K})}
\newcommand{\RepastrQp}{\Rep^{\ast,r}_{\Qp}(\Gamma_{K})}
\newcommand{\Repastrfree}{\Rep^{\ast,r}_{\free}(\Gamma_{K})}
\newcommand{\Repastrtors}{\Rep^{\ast,r}_{\tors}(\Gamma_{K})}
\newcommand{\Repastrft}{\Rep^{\ast, r}_{\finitetype}(\Gamma_{K})}
\DeclareMathOperator{\cone}{cone}
\DeclareMathOperator{\Inv}{Inv}
\DeclareMathOperator{\free}{free}
\DeclareMathOperator{\ceseq}{c\acute{e}}
\DeclareMathOperator{\MP}{P_{\varpi}}
\newcommand{\isoto}{\overset{\sim}\to}
\newcommand{\longisoto}{\overset{\sim}\longto}
\newcommand{\Kbar}{\overline K}
\newcommand{\OK}{\sO_{K}}
\newcommand{\ses}[5]{0\longto #1 \overset{#2}\longto #3 \overset{#4}\longto #5\longto 0}
\newcommand{\tikzses}[6]{\begin{tikzcd} #1 \pgfmatrixnextcell 0  \ar[r] \pgfmatrixnextcell #2 \ar[r, "#3"] \pgfmatrixnextcell #4 \ar[r, "#5"] \pgfmatrixnextcell #6 \ar[r] \pgfmatrixnextcell 0 \end{tikzcd}}
\newcommand{\StacksTag}[1]{\href{http://stacks.math.columbia.edu/tag/#1}{\cite[Tag #1]{stacks-project}}}
\newcommand{\longto}{\longrightarrow}
\newcommand{\sA}{\mathscr A}
\newcommand{\sC}{\mathscr C}
\newcommand{\sM}{\mathscr M}
\newcommand{\sO}{\mathscr O}
\newcommand{\sP}{\mathscr P}
\newcommand{\sQ}{\mathscr Q}
\newcommand{\sT}{\mathscr T}
\newcommand{\bZ}{\mathbb Z}

\newcommand{\Q}{\mathbb Q}
\newcommand{\Z}{\mathbb{Z}}
\newcommand{\Zp}{\Z_{p}}
\newcommand{\Qp}{\Q_{p}}

\newcommand{\fM}{\mathfrak M}
\newcommand{\fS}{\mathfrak S}
\newcommand{\fm}{\mathfrak m}
\newcommand{\fc}{\mathfrak c}

\DeclareMathOperator{\Ext}{Ext}
\DeclareMathOperator{\Gal}{Gal}
\DeclareMathOperator{\GL}{GL}
\DeclareMathOperator{\Hom}{Hom}
\DeclareMathOperator{\id}{id}
\DeclareMathOperator{\Res}{Res}

\newcommand{\Gm}{\mathbb{G}_{m}}
\newcommand{\GmK}{\mathbb{G}_{m,K}}
\newcommand{\GmOK}{\mathbb{G}_{m,\OK}}

\newcommand{\hooklongrightarrow}{\lhook\joinrel\longrightarrow}
\newcommand{\oo}{\infty}

\newcommand{\twoardef}[5]{#1 :\: \left\{ \begin{aligned} #2 \quad& \longrightarrow \quad #3 \\ #4 \quad & \longmapsto \quad    #5 \end{aligned} \right. }
\newcommand{\twoar}[5]{ \begin{aligned} #2 \quad& \overset{#1}\longrightarrow \quad #3 \\ #4 \quad & \longmapsto \quad    #5 \end{aligned} }

\DeclareMathOperator{\fr}{fr}
\DeclareMathOperator{\tors}{tors}
\DeclareMathOperator{\rk}{rk}
\DeclareMathOperator{\pdiv}{pdiv}
\DeclareMathOperator{\DD}{DD}
\DeclareMathOperator{\DEG}{DEG}
\DeclareMathOperator{\stbl}{st}
\DeclareMathOperator{\Crys}{Crys}
\DeclareMathOperator{\fppf}{fppf}
\DeclareMathOperator{\fin}{fin}

\DeclareMathOperator{\pr}{pr}

\DeclareMathOperator{\Rep}{Rep}
\DeclareMathOperator{\Mod}{Mod}

\DeclareMathOperator{\an}{an}

\DeclareMathOperator{\Sh}{Sh}
\DeclareMathOperator{\Sch}{Sch}
\DeclareMathOperator{\colim}{colim}
\DeclareMathOperator{\coker}{coker}

\DeclareMathOperator{\im}{im}

\DeclareMathOperator{\Spec}{Spec}

\DeclareMathOperator{\Ab}{Ab}
\DeclareMathOperator{\cts}{cont}
\DeclareMathOperator{\crys}{crys}
\DeclareMathOperator{\et}{\acute{e}t}
\DeclareMathOperator{\dR}{dR}
\DeclareMathOperator{\fpart}{f}

\begin{document}

\title[Semistable abelian varieties and $\text{Crys}_{1}$]{Semistable abelian varieties and maximal torsion 1-crystalline submodules}
\author{Cody Gunton}
\address{Cody Gunton \\
Department of Mathematical Sciences \\ 
Universitetsparken 5 \\
University of Copenhagen \\
2100 Copenhagen \O \\
Denmark}
\email{cody@math.ku.dk}
\thanks{Supported by the Danish National Research Foundation through the Centre for Symmetry and Deformation (DNRF92).}


\begin{abstract}
Let $p$ be a prime, let $K$ be a discretely valued extension of $\Qp$, and let $A_{K}$ be an abelian $K$-variety with semistable reduction. Extending work by Kim and Marshall from the case where $p>2$ and $K/\Qp$ is unramified, we prove an $l=p$ complement of a Galois cohomological formula of Grothendieck for the $l$-primary part of the N\'eron component group of $A_{K}$. Our proof involves constructing, for each $m\in \Z_{\geq 0}$, a finite flat $\OK$-group scheme with generic fiber equal to the maximal 1-crystalline submodule of $A_{K}[p^{m}]$. As a corollary, we have a new proof of the Coleman-Iovita monodromy criterion for good reduction of abelian $K$-varieties.
\end{abstract}



\maketitle

\normalem

\numberwithin{equation}{section}

\section{Introduction}
\label{sec:introduction}

Let $p$ be any rational prime, let $K/\Qp$ be a discretely valued extension of $\Qp$, let $A_{K}$ be an abelian $K$-variety, and write $A$ for the N\'eron model of $A_{K}$. Let $\Phi_{A_{K}}$ be the finite abelian group of connected components of the special fiber of $A$. Fix an embedding $K\subset \Kbar$ into an algebraic closure and let $\Gamma_{K} = \Gal(\Kbar/K)$. Assume that $A_{K}$ has semistable reduction. We prove a formula for the $p$-primary part $\Phi_{A_{K}}[p^{\oo}]$ of $\Phi_{A_{K}}$. The formula has a much older $l\neq p$ analogue due to Grothendieck. He established the following result without any restriction on the ramification of $K/\Qp$.
  
  \begin{proposition}[Corollary of \cite{grothendieck-monodromy} 11.3.8]
    \label{grothendieck-l-not-p}
  Let $l$ be a prime distinct from $p$. Writing $T_{l}A_{K}$ for the $l$-adic Tate module of $A_{K}$, and $I_{K}\subset \Gamma_{K}$ for the inertia subgroup, there is a canonical isomorphism of unramified $\bZ_{l}[\Gamma_{K}]$-modules
  \[
\Phi_{A_{K}}[l^{\oo}] \simeq (H^{1}_{\cts}(I_{K}, T_{l}A_{K}))_{\tors},
\]
where the cohomology group $H^{1}_{\cts}(I_{K}, T_{l}A_{K})$ is formed using continuous cochains, and where $(-)_{\tors}$ gives the maximal torsion subgroup.
\end{proposition}

Loosely speaking, the formula breaks down when $l=p$ because maximal unramified submodules tend to be too small in that case. The following idea is pursued in the work \cite{kim-marshall} of Kim and Marshall: to `repair' Grothendieck's formula in the case $l=p$, one should replace the condition of unramifiedness by another condition from $p$-adic Hodge theory. Specifically, for $r\in \Z_{\geq 0}$, a torsion $\Z_{p}[\Gamma_{K}]$-module is ``$r$-crystalline'' if it arises as a quotient of a $\Gamma_{K}$-stable $\Zp$-lattice in some crystalline $\Qp[\Gamma_{K}]$-module with Hodge-Tate weights in the interval $[0,r]$.  For any $\Zp[\Gamma_{K}]$-module $M$ of finite-type over $\Zp$ and any $r\in \Z_{\geq 0}$, the notion of a maximal $r$-crystalline submodule $\Crys_{r}(M)$ is well defined, and gives rise to left exact functor $\Crys_{r}(-)$. Now, $0$-crystallinity is simply unramifiedness, and the derived functors of $\Crys_{0}$ compute cohomology of $I_{K}$. From this perspective, it is natural to consider using a derived functor of $\Crys_{1}(-)$ as a replacement for inertial cohomology in the $l=p$ setting. We will prove the following result.

\begin{theorem}[\Cref{main-theorem-component-groups}; \cite{kim-marshall}, Section 2 for $p>2$ and $K$ unramified]
  \label{intro-main-theorem-phi}
  Suppose $A_{K}$ has semistable reduction. Let $R^{1}\Crys_{1}(T_{p}A_{K})$ denote the continuous derived functor of $\Crys_{1}$ described in \Cref{sec:crys-r-derived-functors}. Then there is a canonical isomorphism of unramified $\Zp[\Gamma_{K}]$-modules
  \[
    \Phi_{A_{K}}[p^{\oo}] \simeq (R^{1}\Crys_{1}(T_{p}A_{K}))_{\tors}.
  \]
\end{theorem}

When $A_{K}$ has semistable reduction, $p>2$, and $K/\Qp$ is unramified, Kim and Marshall use Fontaine-Laffaille theory to choose a finite flat $\OK$-group $Q_{p^{m}}$ with generic fiber isomorphic to $\Crys_{1}(A_{K}[p^{m}])$. They then use Raynaud's results on prolongations 
to produce an embedding of the identity component $Q_{p^{m}}^{\circ}$ into the N\'eron model $A$, and to deduce that this embedding gives an isomorphism of $Q_{p^{m}}$ with the finite flat $\OK$-subgroup. Following a strategy in \cite{grothendieck-monodromy}, they construct isomorphisms
\begin{equation}
  \label{eq:intro-phi-and-crys-1}
  \Phi_{A_{K}}[p^{m}] \simeq \frac{\Crys_{1}(T_{p}A_{K}\otimes \Z/p^{m}\Z)}{\Crys_{1}(T_{p}A_{K})\otimes \Z/p^{m}\Z}
\qquad (m \geq 1).
\end{equation}

To remove the assumptions on $p$ and $K$, we take a different approach to constructing this isomorphism. Rather than choosing $Q_{p^{m}}$ by an existence result and demonstrating that it has a certain relationship to $A$, we construct a finite flat $\OK$-group ${}^{*}\sQo[p^{m}]$ directly from the N\'eron model $A$ of $A_{K}$. For this, we will use a semiabelian $K$-variety, the generic fiber of the Raynaud extension $\widetilde A$ attached to $A_{K}$. Using the monodromy pairing, we first construct a finite flat $\OK$-group with monodromy, and then define ${}^{*}\sQo[p^{m}]$ as a subobject of this. Our main technical result, \Cref{main-theorem-torsion}, is the construction of canonical isomorphisms
\begin{equation}
  \label{intro-main-torsion-isom}
  \Crys_{1}(A_{K}[p^{m}])\simeq {}^{*}\sQo[p^{m}]_{K} \qquad (m\geq 1).
\end{equation}
We prove this using a recent result of torsion $p$-adic Hodge Theory, the Torsion Full Faithfulness Theorem of Ozeki, which improves on an older result of Breuil (see \Cref{torsion-ff}). Taking the limit over $m$,  one obtains an isomorphism
\begin{equation}
  \label{intro-main-integral-isom}
  \Crys_{1}(T_{p}A_{K}) \simeq T_{p}\widetilde A_{K}.
\end{equation}
With our descriptions of the ``numerator'' and ``denominator'' in \eqref{eq:intro-phi-and-crys-1}, the isomorphism given in \eqref{eq:intro-phi-and-crys-1} is constructed using results on the monodromy pairing. From here,  we may reuse the homological algebra arguments of \cite{kim-marshall} to prove \Cref{intro-main-theorem-phi}.


Having described the strategy, we further motivate our approach and describe its relationship to existing work.

One knows that a finite flat group scheme with generic fiber isomorphic to $\Crys_{1}(A_{K}[p^{m}])$, like our ${}^{*}\sQo[p^{m}]$, exists by the following result, which postdates the work of Kim and Marshall.

\begin{theorem}[Tate; Raynaud; Kisin; W.\! Kim; see \Cref{1-crys-and-geom-attributions}]
  \label{1-crys-things-are-generic-fibers-intro}
  Let $M$ be a finite, torsion $\Z_{p}[\Gamma_{K}]$-module. Then $M$ is the generic fiber of a finite flat $\OK$-group if and only if $M$ is $1$-crystalline.
\end{theorem}

With \Cref{1-crys-things-are-generic-fibers-intro}, the proof of Kim and Marshall extends immediately to the case where the ramification index $e$ of $K/\Qp$ satisfies $e<p-1$ (in particular, it must be that $p>2$). However, this ramification restriction is essential to their approach, which relies on Raynaud's work on prolongations. Conditions weaker than fullness of the generic fiber functor which are valid without any condition on $e$ have been proved by Bondarko \cite{bondarko-ffgss} and Liu \cite{liu-torsion-fontaine}. These results take the following form: there exists a constant $\fc\in \Z_{\geq 0}$ depending only on $e$ such that, for any two finite flat $\OK$-groups $H'$, $H$ and any morphism $f_{K}: H'_{K}\to H_{K}$, there exist a morphism $H'\to H$ inducing $p^{\fc}f_{K}$ on generic fibers. In the cases where $\fc>0$, we were unable to use these results to generalize the work of Kim and Marshall.

Kisin's proof that torsion 1-crystallinity is the same as finite flatness loses contact with the representation $T_{p}A_{K}$. In particular, the proof uses a functor $\fM(-): \RepcrysoneZp\to \Mod^{1,\fr}_{\fS}$, which is not exact. While one can make\footnote{See \cite{madapusi-pera-log-pdivs-and-sst-reps}; \cite{kisin-f-crystals} does consider categories of $\fS$-modules with monodromy, but the monodromy maps are only defined after inverting $p$.} a variant of $\fM(-)$ amenable to the study of lattices in 1-semistable representations by enriching $\Mod^{1,\fr}_{\fS}$ with monodromy operators, the non-exactness of $\fM(-)$ obstructs some naive approaches to proving \Cref{main-theorem-torsion}. Kisin has shown in \cite{kisin-f-crystals} that there exists a $p$-divisible group $G_{n}[p^{\oo}]$ over $\OK$ such that $G_{n}[p^{m}]_{K}$ is isomorphic to $\Crys_{1}(A_{K}[p^{m}])$. However, there is no immediately useful relationship between the $G_{n}[p^{m}]$ for different values of $m$.

Finite flat group schemes with monodromy were defined by Kato, and are essentially the same as log finite flat group schemes; see the introduction to \Cref{sec:ffgs-monodromy-1sstmodules} for more on this point. We prefer to work with monodromy structures, rather than with log structure, because we feel it clarifies the role of the Torsion Full Faithfulness Theorem.

An isomorphism as in \eqref{intro-main-integral-isom} is also given by Coleman and Iovita (\cite{coleman-iovita-1} 4.6), where it is a key ingredient in the proof of their well known criterion for good reduction of abelian $K$-varieties. Our isomorphisms \eqref{intro-main-torsion-isom} give a torsion refinement that seems to be inaccessible via the $p$-adic integration approach of \cite{coleman-iovita-1}. Using our isomorphism in place of theirs, we obtain a proof of the Coleman-Iovita criterion (\Cref{coleman-iovita-breuil}).

On the issue of removing the semistability assumption in \Cref{intro-main-theorem-phi}, we make no progress. In the work of Kim and Marshall, this assumption is needed to relate $\Phi_{A_{K}}[p^{m}]$ to the finite parts of quasi-finite flat group schemes (see \eqref{eq:phi-p-m-ito-fparts}), and in our work it is needed for the construction of a 1-motive. We have not, however, ruled out the possibility that \Cref{intro-main-theorem-phi} holds in general. To make progress on this question, one would have to understand more about 1-crystalline lifts of residual representations in situations of unipotent reduction. In the case of elliptic curves, the paper \cite{MR3324930} seems to be a useful point of departure.

Grothendieck showed that \Cref{grothendieck-l-not-p} actually holds in greater generality, in particular, for discretely valued $K$ of characteristic $p>0$, and it is natural to ask if an $l=p$ analogue holds over such a field. Our work is restricted to $K$ of characteristic zero at least by its use of Ozeki's Torsion Full Faithfulness Theorem. This theorem is proved using $(\varphi, \widehat{G})$-modules, which are objects of the mixed characteristic world. As a replacement for such objects, one might consider the rigid cohomology of Lazda and P\'{a}l, which has already allowed for the proof of a $p$-adic Néron-Ogg-Shafarevich criterion, Theorem 5.74 in \cite{lazda-pal-rigid-coh}.

We hope we have made our debt to Kim and Marshall apparent. The idea of studying derived functors of $\Crys_{1}$ originates in their work \cite{kim-marshall}, as does their definition of these derived functors in the spirit of Jannsen's definition of continuous \'etale cohomology. Marshall's thesis \cite{marshall-thesis} was an invaluable resource to us in our study of these derived functors and the general categories of (e.g., non-finite type torsion) $r$-crystalline $\Gamma_{K}$-modules that arise in their construction.

The present work is derived from the author's PhD thesis, which was written under the guidance of Bryden Cais. The author sincerely thanks him for suggesting the project of removing the ramification restriction in the work of Kim and Marshall, and for his consistent support and help. The author thanks Keerthi Madapusi Pera for sharing his unpublished preprint \cite{madapusi-pera-log-pdivs-and-sst-reps}. He also thanks Lars Halvard Halle for his comments on an earlier draft. Finally, the author thanks the referee for mathematical corrections and suggestions relating to the exposition that have substantially improved the quality of the article.
\subsection{Notation}
\label{subsec:notation}

Let $p$ be any prime and let $K/\Qp$ be a discretely valued field extension. Write $\fm_{K}\subset \OK$ be the maximal ideal, let $\kappa$ be the residue field of $\OK$, and let $v_{K}$ be the valuation of $\OK$ satisfying $v_{K}(K^{\times}) = \bZ$. Fix an embedding $K\subset \Kbar$ into an algebraic closure, let $\Gamma_{K} := \Gal(\Kbar/K)$, and let $C$ be the completion of $\Kbar$ under the absolute value extending the absolute value of $K$.

Fix a uniformizer $\varpi\in \OK$. Let $\sO_{C}^{\flat} := \lim_{x\mapsto x^{p}}\sO_{C}$, let $\varpi^{\flat} = (\varpi^{\flat,(i)})_{i\geq 0}\in \sO_{C}^{\flat}$ be an element with $\varpi^{\flat, (0)}=\varpi$, and let $\epsilon = (\epsilon^{(i)})_{i\geq 0}\in \sO_{C}^{\flat}$ be an element with $\epsilon^{(0)}= 1 \neq \epsilon^{(1)}$. Then $\varpi^{\flat}$ is a system of $p$-power roots of $\varpi$ and $\epsilon$ is a nontrivial system of $p$-power roots of 1. Let $\Gamma_{\varpi^{\flat}}\subset \Gamma_{K}$ be the (non-normal) subgroup consisting of those automorphisms that fix each of the $\varpi^{\flat, (i)}$.

\subsection{Example: Tate curves}
\label{subsec:q-tate-eg}
We describe our results in the special case of a Tate curve. While this very simple example could be treated in a more streamlined and direct way, we give a presentation that aims to motivate our more general constructions. To that end, we use notation and terminology that will not be properly defined until later.

Let $q$ be an element of $K^{\times}$ satisfying $v_{K}(q)>0$, and let $E_{q,K}$ be the associated Tate curve. It is well known that the homomorphism $\iota_{q,K}: \bZ_{K}\longto \GmK$ of $K$-group schemes sending $1$ to $q$ fits into a short exact sequence of $K$-rigid spaces
\begin{align}
  \label{eq:q-tate-unif-ses}
  \ses{\bZ_{K}^{\an}}{\iota_{q,K}}{
  \GmK^{\an}}{\psi_{K}}{E_{q,K}^{\an}},
\end{align}
where $\psi_{K}$ is given by explicit series over $K$. Let $n\geq 2$ be a power of $p$. The short exact sequence \eqref{eq:q-tate-unif-ses}  gives rise to the short exact sequence of finite $K$-group schemes
\begin{align}
  \label{eq:q-tate-mot-ses}
 \ses{\mu_{n,K}}{\psi_{K}}{E_{q,K}[n]}{}{\Z_{K}/n\Z_{K}}.
\end{align}
Since the $p$-adic cyclotomic character is 1-crystalline, $\Crys_{1}(E_{q,K}[n])$, the maximal 1-crystalline subobject of $E_{q, K}[n]$, is at least as large as $\Z/n\Z(1)$. The map $\iota_{q,K}$ can be regarded a log 1-motive over $K$
\begin{align}
  \label{eq:q-tate-1-motive}
  \sT_{q,K} = [\Z_{K}\overset{\iota_{q,K}}\longto \GmK],
\end{align}
which is, in particular, a complex of $K$-groups concentrated in degrees $-1$ and $0$. The integer $n$ gives a morphism of complexes $[n]: \sT_{q,K}\to \sT_{q,K}$. Defining $\sT_{q,K}[n] := H^{-1}(\cone([n]))$, $\psi_{K}$ induces an isomorphism (cf.~\eqref{eq:1-mot-n-tors-isom})
\begin{align*}
  \sT_{q,K}[n]\simeq E_{q,K}[n].
\end{align*}
As a torsion $\Gamma_{K}$-module, $E_{q,K}[n]$ is 1-semistable. To describe the maximal 1-crystalline subobject $\Crys_{1}(E_{q,K}[n])$, we will use Raynaud's decomposition of  $\sT_{q,K}$ (see \Cref{raynaud-decomp}). Use the uniformizer $\varpi\in \OK$ to write $q = u\varpi^{v}$ with $u\in \OK^{\times}$ and $v\in \Z_{\geq 1}$. This allows us to define two log 1-motives
  \begin{alignat*}{2}
      \sT_{q,K,\varpi}^{1}
      &:= [\Z_{K}\overset{\iota_{q,K,\varpi}^{1}}\longto \GmK],
      &&\qquad \iota_{q,K,\varpi}^{1}(1) = u
      \\
      \sT_{q,K,\varpi}^{2}
      &:= [\Z_{K}\overset{\iota_{q,K,\varpi}^{2}}\longto \GmK],
      &&\qquad \iota_{q,K,\varpi}^{2}(1) = \varpi^{v}.
  \end{alignat*}
Then $\iota_{q,K} = \iota_{q, K, \varpi}^{1} \cdot \iota_{q, K, \varpi}^{2}$ using the group law in $\GmK$. It is clear that $\sT_{q,K,\varpi}^{1}$ extends to an $\OK$-1-motive $\sT_{q,\varpi}^{1}$. Defining torsion in 1-motives as we did above, we have a finite flat $\OK$-group scheme
      $\sT_{q,\varpi}^{1}[n]$ with generic fiber $\sT_{q,K,\varpi}^{1}[n]$, so that the $\Gamma_{K}$-module $\sT_{q,K,\varpi}^{1}[n]$ is 1-crystalline.

     Under the isomorphism $c_{(-),n}: K^{\times}/(K^{\times})^{n}\to H^{1}(\Gamma_{K}, \Z/n\Z(1))$ of Kummer Theory, the class $\eta(\sT_{q,K}, n)$ (resp., $\eta(\sT_{q,K,\varpi}^{1}, n)$, resp., $\eta(\sT_{q,K,\varpi}^{2}, n)$) of $\sT_{q,K}$ (resp., $\sT_{q,K,\varpi}^{1}$, resp., $\sT_{q,K,\varpi}^{2}$) in $\Ext^{1}_{\ModZpGK}(\Z/n\Z, \Z/n\Z(1))$  arises from $q$ (resp., $u$, resp., $\varpi^{v}$). The equation $c_{q,n} = c_{u,n} + c_{\varpi^{v},n}$ gives a Baer sum decomposition
\begin{align}
\label{eq:q-tate-tors-dec}
  \eta(\sT_{q,K}, n)
  = \eta(\sT_{q,K,\varpi}^{1}, n)
  + \eta(\sT_{q,K,\varpi}^{2}, n).
\end{align}
By construction of the Baer sum, $\sT_{q,K}[n]$ is a subquotient of $\sT_{q,K,\varpi}^{1}[n] \oplus \sT_{q,K,\varpi}^{2}[n]$. Since finite direct sums and subquotients of 1-crystalline $\Gamma_{K}$-modules are again $1$-crystalline, taking a difference in \eqref{eq:q-tate-tors-dec} allows us to conclude that $E_{q,K}[n]$ is 1-crystalline if and only if $\sT_{q,K,\varpi}^{2}[n]$ is 1-crystalline.

Since $\sT_{q,K,\varpi}^{2}[n]$ corresponds to the 1-cocycle $c_{\varpi^{v},n} = v c_{\varpi,n}:\Gamma_{K}\to \Z/n\Z(1)$, we see that $\sT_{q,K,\varpi}^{2}[n]$ is split, i.e., $\sT_{q,K,\varpi}^{2}[n] \simeq \Z/n\Z(1)\oplus \Z/n\Z$ as a $\Gamma_{K}$-module, if $n\mid v$. In particular, in this case $\sT_{q,K,\varpi}^{2}[n]$ is 1-crystalline. Conversely, by the Torsion Full Faithfulness Theorem of Ozeki, Breuil and others (see \Cref{torsion-ff}), if $\sT_{q,K,\varpi}^{2}[n]$ is 1-crystalline and $V$ is another 1-crystalline $\Gamma_{K}$-module, then any isomorphism \[\Resoo(\sT_{q,K,\varpi}^{2}[n]) \simeq \Resoo V\]
induces an isomorphism $\sT_{q,K,\varpi}^{2}[n]\simeq V$. But $c_{\varpi, n}(\Gamma_{\varpi^{\flat}}) = 0$, so any choice of an adapted basis for $\sT_{q,K,\varpi}^{2}[n]$ induces an isomorphism
\[\Resoo \sT_{q,K,\varpi}^{2} [n]
  \simeq \Resoo
  \Big(\bZ/n\bZ(1)\oplus \bZ/n\bZ\Big).
\]
We conclude that $E_{q,K}[n]$ is 1-crystalline if and only if $v_{p}(n) \leq v_{p}(v_{K}(q))$.

We push this argument to describe $\Crys_{1}(E_{q,K}[n])$ in general. Suppose $v_{p}(n)\geq v_{p}(v_{K}(q))$. With reference to \ref{eq:q-tate-mot-ses}, form a basis of $E_{q,K}[n]$ by choosing an element $x_{1}$ that spans the image of $\psi_{K}$ and letting $x_{2}\in E_{q,K}[n]$ be a lift of $1\in \bZ/n\bZ$. Then $\Crys_{1}(E_{q,K}[n])$ is generated by $x_{1}$ and $p^{w}x_{2}$ for some integer $w$ satisfying $0 \leq w \leq v_{p}(n)$. In particular, $p^{w}E_{q,K}[n]$ is contained in $\Crys_{1}(E_{q,K}[n])$. Since 1-crystallinity is inherited by submodules, we see that $E_{q,K}[p^{v_{p}(n)-w}]$ is 1-crystalline. Our above work shows that $v_{p}(n)-w\leq v_{p}(v_{K}(q))$, and the maximality property of $\Crys_{1}(E_{q,K}[n])$ implies that $w = v_{p}(n)-v_{p}(v_{K}(q))$. Therefore, for $v_{p}(n) \geq v_{K}(q)$, we have $\Crys_{1}(E_{q,K}[n]) = \bZ/n\bZ x_{1} + p^{v_{p}(n)-v_{p}(v_{K}(q))}\bZ/n\bZ x_{2}$.

The subspace $\im \psi_{K}\simeq \bZ/n\bZ(1)$ is distinguished by the fact that it arises from a 1-crystalline submodule of $T_{p}E_{q,K}$. A direct calculation with the period ring $B_{\crys}$ show that $T_{p}E_{q,K}$ is not crystalline. Being that $T_{p}E_{q,K}$ is of rank 2 over $\Zp$, we see that $\Crys_{1}(T_{p}E_{q,K})$ is equal to the copy $\Zp(1)$ obtained in the limit from the maps $\psi_{K}$ for varying $n$. It is well known that the N\'eron component group $\Phi_{E_{q}}$ attached to $E_{q,K}$ is cyclic of order $v_{K}(q)$, hence that the $p$-Sylow subgroup $\Phi_{E_{q}}[p^{\oo}]$ is cyclic of order $p^{v_{p}(v_{K}(q))}$. Altogether, we have the following result.

\begin{proposition}
  \label{tate-curve-prop}
  Let $q\in K^{\times}$ satisfy $v_{K}(q) > 0$, let $m\in \bZ_{\geq 1}$, and let $v = v_{p}(\# \Phi_{E_{q}}[p^{\oo}])$. Then
  \begin{equation}
    \label{eq:q-tate-crys-1-cases}
    \Crys_{1}(E_{q,K}[p^{m}]) = 
    \begin{cases}
      E_{q,K}[p^{m}] & \text{ if } m\leq v \\
      \bZ/p^{m}\bZ(1) + p^{m-v}E_{q,K}[p^{m}]
      & \text{ if } m > v
    \end{cases}
  \end{equation}
\end{proposition}
\noindent We will see later that there is an isomorphism of finite abelian groups
  \begin{align*}
        \frac{
      \Crys_{1}(T_{p}E_{q,K}\otimes_{\Zp} \bZ/p^{m}\bZ)
    }{\Crys_{1}(T_{p}E_{q,K})\otimes_{\Zp} \bZ/p^{m}\bZ}
    \simeq \Phi_{E_{q}}[p^{m}].
  \end{align*}
  We end this example by noting that the cases in \eqref{eq:q-tate-crys-1-cases} can be described uniformly. Let $N_{p^{m}}^{*}$ be the endomorphism of the finite flat $\OK$-group scheme $\bZ/p^{m}\bZ$ given by multiplication by $v_{K}(q)$. Let ${}^{*}\sT^{1}_{\varpi}[p^{m}]$ denote
the pullback of the class of $\sT^{1}_{\varpi}[p^{m}]$ in $\Ext(\bZ/p^{m}\bZ, \bZ/p^{m}\bZ(1))$ by the inclusion $\ker N_{p^{m}}^{*}\subset \bZ/p^{m}\bZ$. Then, for all $m\in \bZ_{\geq 1}$, there is an isomorphism
\[{}^{*}\sT^{1}_{\varpi}[p^{m}]_{K}\simeq \Crys_{1}(E_{q,K}[p^{m}]).\]

\section{Finite flat group schemes with monodromy and log 1-motives}
\label{sec:monodromy}

Throughout this section we adopt the notation of \Cref{subsec:notation}.

We recall some known results about semistable degenerations of abelian varieties; nothing is original beyond our presentation. Our main goals are to show that the torsion in a semistable abelian $K$-variety is the torsion in a certain log 1-motive; to describe Raynaud decomposition of log 1-motives; and to record the relationship between the monodromy pairing and the N\'eron component group.

\subsection{Degenerations of abelian varieties}

Let $A_{K}$ denote a semistable abelian $K$-variety. The N\'eron model $A$ of $A_{K}$ is a smooth $\OK$-group. Its special fiber $A_{\kappa}$ need not be connected. We are concerned with $A_{\kappa}/A_{\kappa}^{\circ}$, which (\cite{milne-alg-gps-CUP}, 1.1.33) is a finite \'etale $\kappa$-group. We define two $\OK$-groups.

\begin{definition}
  \label{def:connected-neron}
Define $A^{\circ}$ to be the open subgroup scheme $A^{\circ}\subset A$ with generic fiber $A_{K}$ and with special fiber equal to the identity component of $A_{\kappa}$.

  Define the N\'eron component group of $A_{K}$ to be the unique \'etale $\OK$-group $\Phi_{A_{K}}$ with special fiber $\Phi_{A_{K},\kappa} = A_{\kappa}/A_{\kappa}^{\circ}$.
\end{definition}

Since $\kappa$ is perfect, the identity component $A_{\kappa}^{\circ}$ lies in a canonical short exact sequence (\cite{milne-alg-gps-CUP} 8.27 and 16.15)
\begin{align}
  \label{chevalley-seq-kappa}
  \ses{T_{\kappa}\times_{\kappa} U_{\kappa}}{}{A^{\circ}_{\kappa}}{}{B_{\kappa}}
\end{align}
where $B_{\kappa}$, $U_{\kappa}$ and $T_{\kappa}$ are algebraic $\kappa$-group which are, respectively, abelian, unipotent, and a torus. We say that $A_{K}$ has semistable reduction (or simply that $A_{K}$ is semistable) if $U_{\kappa}$ is trivial. By \cite{grothendieck-monodromy} 2.2.1, $A_{K}$ has semistable reduction if and only if, for all $n\in \Z_{\geq 1}$, the multiplication-by-$n$ map on $A^{\circ}$ is quasi-finite, flat, and surjective. When these conditions are satisfied, the natural injection
  \[
    A_{\kappa}[n]/A_{\kappa}^{\circ}[n]
             \longto \Phi_{A_{K},\kappa}[n]
  \]
  is surjective, opening up the possibility of studying the N\'eron component group (a finite abelian group with unramified $\Gamma_{K}$-action) via the quasi-finite flat $\OK$-group $A[n]$. Rather than to work directly with $A$, we will work with $A^{\circ}$, the Mumford degeneration data formed from this, and Grothendieck's monodromy pairing.

  The theory of degenerations developed by Mumford and Faltings-Chai is exposed in great detail in Lan's book \cite{lan-PUP}. We give only a summary of background we will need, omitting many technical issues relating to, for example, independence of certain choices that are made in order to define the functors $F$ and $M$ below. We note, as Lan does, that most of the results we cite can be found in the older works \cite{raynaud-uniformization-icm}, \cite{mumford-degen-abelian-varieties}, \cite{faltings-chai} and \cite{bosch-lutkebohmert-degen}.  

  The Raynaud extension attached to $A_{K}$ is obtained by successively lifting \eqref{chevalley-seq-kappa} (in the case $U_{\kappa}=0$) to the rings $\sO_{K}/\fm_{K}^{i+1}$ or all $i\geq 0$, and applying Grothendieck's theorem that the resulting semiabelian formal $\OK$-group scheme is algebraizable (see \cite{lan-PUP}, 3.3.3.6 and 3.3.3.9) to a group $\widetilde A$, the Raynaud extension\footnote{Letting $G = A^{\circ}$, Lan uses the notation $G^{\natural}$. This is what Grothendieck would call $G^{\natural\circ}$, reserving $G^{\natural}$ for an associated $G^{\natural\circ}$-torsor over the component group $\Phi_{A_{K}}$. To avoid confusion, since we will study the component group, we use the notation of Faltings-Chai, though we adopt the abuse of writing $\widetilde A$ for what should be called $\widetilde{A^{\circ}}$.}. As just described, this fits into a short exact sequence
\begin{align}
  \label{tilde-A-ses}
  \ses{T}{}{\widetilde A}{}{B}
\end{align}
with $B$ an abelian $\OK$-scheme and $T$ an $\OK$-torus. Let $Y^{D}$ be the character group of $T$. Replacing $A_{K}$ by its dual in the above discussion gives the dual Raynaud extension $\widetilde A^{D}$, which is an extension of the dual $B^{D}$ of $B$ by a torus $T^{D}$. Let $Y$ be the character group of $T^{D}$, a twisted constant group which becomes free of finite rank after finite \'etale base change on $\OK$. From the rigid analytic perspective, $A_{K}$ can be recovered as the cokernel of the analytification of a homomorphism of $K$-group schemes $Y_{K}\to \widetilde A_{K} := (\widetilde A)_{K}$. An important, concrete consequence (known in a different guise since \cite{grothendieck-monodromy}, at the latest) is that $A_{K}[n]$ gives a class in $\Ext(Y_{K}/nY_{K}, \widetilde A_{K}[n])$. It was shown in \cite{mumford-degen-abelian-varieties} (in the case $B_{\kappa}=0$) and \cite{faltings-chai} (in the general case) that a construction of $A_{K}$ from $Y_{K}\to \widetilde A_{K}$ and of a resulting extension structure on $A_{K}[n]$ can be given in algebraic terms which are valid for degenerations over rings more general than $\OK$. For this, we introduce a category of degeneration data.

By \cite{lan-PUP} 3.1.5.1, the $T$-torsor (resp., $T^{D}$-torsor) $\widetilde A$ (resp., $\widetilde A^{D}$) is the same information as a homomorphism $c^{D}:Y \to B$ (resp., $c: Y^{D}\to B^{D}$). The pullback of the Poincar\'e sheaf $\sP_{B}^{-1}$ on $B\times_{\OK} B^{D}$ to $Y\times_{\OK}Y^{D}$ under $c^{D} \times c$ need not admit a trivialization, but its restriction to the generic fiber must.

\begin{definition}[\cite{lan-PUP} 4.4.1 and 4.4.10]
  \label{DEG-DD-def}
  Let $\DEG(\OK,\fm_{K})$ denote the category of $\OK$-schemes whose objects are the semiabelian $\OK$-schemes $G$ such that $G_{K}$ is an abelian variety.

  Let $\DD(\OK,\fm_{K})$ denote the category of triples $(B, Y, Y^{D}, c, c^{D}, \tau)$, such that
  \begin{enumerate}
  \item $B$ is a semiabelian $\OK$-scheme;
  \item $Y$ and $Y^{D}$ are twisted constant $\OK$-groups  which are isomorphic to the constant group $\Z^{t}$ after a finite \'etale base change on $\OK$, for some $t\in \Z_{\geq 0}$;
  \item $c^{D}:Y\to B$ and $c: Y^{D}\to B^{D}$ are homomorphisms; and
  \item $\tau$ is a trivialization of the $\GmK$-biextension $(c^{D}\times c)^{*}\sP_{B,K}^{-1}$;
  \end{enumerate}
  where $\tau$ satisfies the following positivity condition:
  \begin{itemize}
  \item[$(\ast)$] there exists an injection with finite cokernel $\phi:Y\to Y^{D}$ such that, for all $y\in Y_{K}$, the section $\tau(y,\phi(y))\in (c^{D}(y)\times c(\phi(y)))^{*}\sP_{B,K}^{\otimes -1}$ extends to a section of $(c^{D}(y)\times c(\phi(y)))^{*}\sP_{B}^{\otimes -1}$, and for $y\neq 0$, the morphism $(c^{D}(y)\times c(\phi(y)))^{*}\sP_{B}^{\otimes}\to \sO_{\Spec \OK}$ induced by $\tau(y, \chi)$ factors through $\fm_{K}$.
  \end{itemize}

By \cite{lan-PUP} 4.2.1.7, such a $\tau$ determines and is determined by a period homomorphism, i.e., a homomorphism $\iota_{K}:Y_{K}\to \widetilde A_{K}$ lifting $c^{D}:Y_{K}\to B_{K}$ under the surjection $\widetilde A\to B$ of \eqref{tilde-A-ses}.
\end{definition}

Let $\tau$ be as in \Cref{DEG-DD-def}. Indexing the components of $Y_{K}\times_{K}Y^{D}_{K}$ by $K$-points $(y,\chi)$, the datum of $\tau$ is equivalent to a collection of sections $\tau(y,\chi) \in (c^{D}(y)\times c(\chi))^{*}\sP_{B, K}^{-1}$ over the components of $Y_{K}\times_{K}Y^{D}_{K}$. The image of
\[
  (c^{D}(y)\times c(\chi))^{*} \sP_{B}
  \subset (c^{D}(y)\times c(\chi))^{*}\sP_{B, K}
  \overset{\tau(y,\chi)}\longisoto \sO_{\Spec K}
\]
is the same information as an $\OK$-submodule
\begin{align}
  \label{eq:I-y-chi-def}
I_{y,\chi}\subset K  
\end{align}
Writing $I_{y,\chi} = \fm_{K}^{\mu(y,\chi)}$
defines a canonical map $\mu: Y(\Kbar) \times  Y^{D}(\Kbar)\longto \Z$ which is $\Gamma_{K}$-equivariant. Since $\tau$ is a trivialization of $\sP_{B,K}^{-1}$ as a $\GmK$-biextension, one obtains a bilinear morphism
\begin{equation}
  \label{eq:mu-def}
  \mu: Y \times_{\OK} X\longto \Z_{\OK}.
\end{equation}

Using enrichments of the categories $\DEG(\OK, \fm_{K})$ and $\DD(\OK, \fm_{K})$, it can be shown (\cite{lan-PUP} 4.5.5.5) that there are inverse equivalences of categories \[F(\OK,\fm_{K}):
  \DEG(\OK, \fm_{K})\longto \DD(\OK, \fm_{K})\]
and
\[M(\OK,\fm_{K}):
  \DD(\OK, \fm_{K})\longto \DEG(\OK, \fm_{K}).\]

In order to apply results of \cite{raynaud-1-motifs} and \cite{bertapelle-candilera-cristante}, we  will view $F(\OK,\fm_{K})(A^{\circ})$ as a log 1-motive. We introduce this notion now.

\begin{definition}
  Let $S$ be a scheme. An \textit{$S$-1-motive} is a complex $\sM = [\iota:Y\to C]$ of commutative $S$-group schemes in degrees $-1, 0$ such that
  \begin{enumerate}[(i)]
  \item $Y$ is a twisted constant group which is \'etale-locally modeled on a free abelian group of finite rank;
  \item $C$ is an extension of an abelian $S$-scheme by an $S$-torus;
  \item $\iota$ is an homomorphism of $S$-groups.
  \end{enumerate}

Suppose now $S = \Spec \OK$. A \textit{log 1-motive} $\sQ$ over $\OK$ is a triple $(Y, C,\iota_{K})$ where $Y$ and $C$ are commutative group schemes over $S$ satisfying (i) and (ii) above and $\iota_{K}$ is a homomorphism $\iota_{K}:Y_{K}\to C_{K}$. In this case $\sQ_{K} := [\iota_{K}:Y_{K}\to C_{K}]$ is a $K$-1-motive. 
\end{definition}
\begin{remark}
  The sort of semiabelian $S$-schemes that can give a term in a 1-motive over $S$ have constant toric rank. In particular, if $A$ is the N\'eron model of an abelian $K$-variety $A_{K}$ with bad reduction, then $A^{\circ}$ is not a term in any 1-motive over $\OK$.
\end{remark}

\begin{remark}
  \label{logarithmic-enhancements}
  Let $S^{\log}$ be the log scheme obtained by equipping $\Spec \OK$ with its canonical (divisorial) log structure. The definition of a log 1-motive does not involve any log schemes. However, from a log 1-motive it is possible to build a complex of logarithmic group schemes defined over $S^{\log}$ by taking a ``logarithmic enhancement'' of the semiabelian scheme $C$. For more on this, see \cite{kkn-4} 1.4.  
\end{remark}

\begin{example}
  A fundamental example of a log 1-motive is given in \eqref{eq:q-tate-1-motive}.
\end{example}

We are interested in the ``torsion'' in the $K$-1-motives underlying log 1-motives, which is defined in the following proposition.

 \begin{proposition}[\cite{andreatta-baerbieri-viale-1-motives}, \S1.3]
  \label{torsion-in-1-motive}
  Let $S$ be a scheme. Let $\sM = [\iota:Y\to C]$ be a 1-motive over $S$. Each positive integer $n$ gives a morphism of complexes of $S$-groups $[n]:\sM\to \sM$ which is multiplication by $n$ on each term. Let $c(\sM, n)$ be the cone of this morphism, and let $\sM[n]:= H^{-1}(c(\sM, n))$. Then $\sM[n]$ is a finite locally free $S$-group scheme, and $(\sM[p^{m}])_{m\geq 1}$ forms a $p$-divisible group $\sM[p^{\oo}]$.
\end{proposition}

Continuing the notation of the proposition, in local sections we have
\begin{align}
  \label{eq:1-mot-n-tors-isom}
    \sM[n] = \frac{\{
    (g,y)\in C \times Y:\: ng = -\iota(y)
    \}}{\{
    (-\iota(y), ny) :\: y\in Y
  \}}.
\end{align}
By divisibility of semiabelian schemes (\cite{blr}, page 180),
the natural homomorphism $\sM[n] \to Y/nY$ is surjective. This fits into a canonical short exact sequence
  \begin{align}
    \label{eq:eta-M-n}
    \tikzses{\eta(\sM, n):}{C[n]}{}{\sM[n]}{}{Y/nY}.
  \end{align}

  \begin{proposition}
    \label{1-motive-gives-torsion}
    Continue notation, so that $A_{K}$ is an abelian $K$-variety with semistable reduction, and $A$ is its N\'eron model. Let
    \[
      \sQ_{K}:= [Y_{K}\overset{\iota_{K}}\to \widetilde A_{K}]
    \]
    be the log 1-motive attached to $A_{K}$ via the functor $F(\OK, \fm_{K})$.
    There is a canonical isomorphism $\sQ_{K}[n] \isoto A_{K}[n]$. Therefore, the short exact sequence $\eta(\sQ_{K},n)$ gives rise to canonical short exact sequence
  \begin{align*}
    \tikzses{\eta(A_{K}, n):}{\widetilde A_{K}[n]}{}{A_{K}[n]}{}{Y_{K}/nY_{K}}.
  \end{align*}
    \proof{
      In this case, \eqref{eq:1-mot-n-tors-isom} reads
  \begin{align*}
    \sQ_{K}[n] = \frac{\{
    (g,y)\in \widetilde A_{K} \times Y_{K}:\: ng = -\iota_{K}(y)
    \}}{\{
    (-\iota_{K}(y), ny) :\: y\in Y_{K}
  \}}.
\end{align*}
Comparing this with $\coker(\iota_{K}^{\an})[n]$, it is clear that rigid analytic uniformization gives a canonical isomorphism $\sQ_{K}[n] \simeq A_{K}[n]$. \qed
    }
  \end{proposition}

  \begin{lemma}
    \label{period-homs-are-injective}
    With notation as in \Cref{1-motive-gives-torsion}, $\ker \iota_{K} = 0$.
    \proof{
      We treat the case where $Y$ is split. The general case follows by \'etale descent. Let $y\in Y(\OK)$. Using $\mid-\mid$ to represent total spaces of invertible sheaves, for each $\chi\in Y^{D}(\OK)$, define a line bundle by pullback as in
\begin{equation*}
  \begin{tikzcd}
    \mid\sO_{\chi}(c^{D}(y))\ \mid \ar[d]\ar[r] &
    \mid\sO_{\chi}\mid \ar[d]\ar[r] &
    \mid\sP_{B}\mid \ar[d]
    \\
    \Spec \OK \ar[r, "(c^{D}(y){, } \id)"] &
    B \times_{\OK} \Spec \OK \ar[r,"\id\times c(\chi)"] &
    B \times_{\OK}B^{D}
  \end{tikzcd}
\end{equation*}

Then $\sO_{\widetilde A} = \bigoplus_{\chi}\sO_{\chi}$ as an $\sO_{B}$-algebra by \cite{lan-PUP} 4.5.2.3. By the proof of \cite{lan-PUP} 3.1.4.4, the identity section of $\widetilde A$ is characterized by the condition that, for each $\chi$, it restricts to the generic fiber of the morphism $e_{\chi,K}^{*}:\sO_{\chi}\to e_{B}^{*}\sO_{\chi}\simeq \sO_{\Spec \OK}$, where the isomorphism is the rigidification arising from the natural rigidification of $\sP_{B}$ along the map $B\to B\times_{\OK} B^{D}$. Each $\tau(y,\chi)\in \sO_{\chi}(c^{D}(y))^{\otimes -1}$ is a morphism $\sO_{\chi}(c^{D}(y))\to \sO_{\Spec K}$. Using the natural maps $\sO_{\chi}\to \sO_{\chi}(c^{D}(y))$ of the above diagram, the map $\iota_{K}(y)$ is defined to be the composition
\[
 \bigoplus_{\chi}\sO_{\chi}
 \longto \bigoplus_{\chi}\sO_{\chi}(c^{D}(y))
 \overset{\sum_{\chi}\tau(y, \chi)}\longto \sO_{\Spec K}.
\]
Therefore, if $y\in \ker\iota_{K}$, then for all $\chi$ we have a commuting diagram
\begin{equation*}
  \begin{tikzcd}
    \sO_{\chi, \Spec K} \ar[r]
    \ar[rr, bend left=20, "\tau(y{,} \chi)"] &
    e_{B}^{*}\sO_{\chi, \Spec K} \ar[r,"\sim"] &
    \sO_{\Spec K} 
    \\
    \sO_{\chi} \ar[r] \ar[u] &
    e_{B}^{*}\sO_{\chi} \ar[r,"\sim"] \ar[u] &
    \sO_{\Spec \OK} \ar[u]
  \end{tikzcd}
\end{equation*}
Since the diagram commutes, the $\OK$-module $I_{y,\chi}$ defined in \eqref{eq:I-y-chi-def} satisfies $I_{y,\chi} = \OK$. In particular, $\tau(y,\chi)$ does not factor through the ideal sheaf $\fm_{K}$. Taking $\phi$ as in \Cref{DEG-DD-def} and setting $\chi = \phi(y)$, the positivity condition of the definition implies that $y = 0$.
      \qed}
  \end{lemma}

  \begin{example}
    In the case where $A_{K} = E_{q,K}$ is the $q$-Tate curve over $K$, the log 1-motive $\sQ_{K} = \sT_{q,K}$ is defined in \eqref{eq:q-tate-1-motive} and $\eta(E_{q, K},n)$ coincides with the short exact sequence \eqref{eq:q-tate-mot-ses}.
  \end{example}

  \subsection{Raynaud decomposition and monodromy}
  
  As we shall see, the exact sequence $\eta(A_{K}, n)$ exhibits $A_{K}[n]$ as an extension of a $0$-crystalline representation by a $1$-crystalline representation, hence as an extension of generic fibers of finite flat $\OK$-group schemes. If $A_{K}$ has bad reduction, then Coleman-Iovita's monodromy criterion for good reduction (\Cref{coleman-iovita-breuil}) coupled with work of Liu (\Cref{liu-breuil-conjecture}) shows that, for $n\gg 0$ a power of $p$, $A_{K}[n]$ is not the generic fiber of a finite flat $\OK$-group. This is related to the fact that, though the individual terms in $\sQ_{K}$ do extend to group schemes over $\OK$, the morphism $\iota_{K}$ does not, in general, extend to a morphism of $\OK$-group schemes. As noted in \cite{raynaud-1-motifs} 2.4.1, the failure of $\iota_{K}$ to extend over $\OK$ is due to nontriviality of the pullback $(c^{D}\times c)^{*}\sP_{B}$ as a $\Gm$-torsor on the $\OK$-scheme $Y\times_{\OK}Y^{D}$. This nontriviality is, in turn, captured by the nontriviality of the monodromy pairing $\mu$. Raynaud's decomposition theorem trivializes the monodromy pairing to get an $\OK$-1-motive, and records the lost information in a $K$-1-motive.

  \begin{proposition}[\cite{raynaud-1-motifs} 4.5.1; see also \cite{bertapelle-candilera-cristante} Theorem 9 and Lemma 12]
  \label{raynaud-decomp}
  Let $\sQ_{K} = F(\OK,\fm_{K})(A^{\circ}) = [Y_{K}\overset{\iota_{K}}\longto \widetilde A_{K}]$ be as in \Cref{1-motive-gives-torsion}. Let $\varpi$ be a uniformizer of $\OK$. Then there are homomorphisms $\iota_{K,\varpi}^{1},\,\iota_{K,\varpi}^{2}:Y_{K}\to \widetilde A_{K}$ such that, using addition in the common target $\widetilde A_{K}$, we have
  \begin{align*}
    \iota_{K} = \iota_{K,\varpi}^{1} + \iota_{K,\varpi}^{2}.
  \end{align*}
  These homomorphisms satisfy
  \begin{enumerate}[(i)]
  \item $\iota_{K,\varpi}^{1}$ is the generic fiber of an $\OK$-1-motive $\sQo:= [Y\overset{\iota^{1}_{\varpi}}\longto \widetilde A]$; and
  \item $\iota_{K,\varpi}^{2}$ factors through the inclusion $ T_{K}\to \widetilde A_{K}$.
  \end{enumerate}
  Let $\sQ_{K,\varpi}^{1}$ and $\sQ_{K,\varpi}^{2}$ be the $K$-1-motives attached to $\iota_{K,\varpi}^{1}$ and $\iota_{K,\varpi}^{2}$, respectively. Then we have an equality in $\Ext(Y_{K}/nY_{K}, \widetilde A_{K}[n])$
  \begin{equation}
    \label{eq:raynaud-baer-decomp}
 \eta(\sQ_{K}, n) = \eta(\sQ_{K,\varpi}^{1}, n) + \eta(\sQ_{K,\varpi}^{2}, n).
  \end{equation}
\end{proposition}

\begin{example}
  In the case of a Tate curve, this is the decomposition of \eqref{eq:q-tate-tors-dec}.
\end{example}

Since $\eta(\sQ_{K,\varpi}^{1}, n)$ is the generic fiber of a sequence of finite flat $\OK$-group schemes $\eta(\sQ_{K}^{1}, n)$, roughly speaking, \eqref{eq:raynaud-baer-decomp} describes how to isolate the non-1-crystallinity of $A_{K}[n]$ in the term $\eta(\sQ_{K,\varpi}^{2},n)$. In the next section, we will develop some of Kato's theory of finite flat group scheme with monodromy, which develops this perspective. We now recall the relationship between the monodromy pairing and the component group (valid only under the semistability assumption) first articulated in \cite{grothendieck-monodromy}.

\begin{definition}
  \label{def:finpR}
Let $R$ be a field or a Henselian discrete valuation ring. Let $\fin^{p}_{R}$ denote the category of commutative $\Spec R$-group schemes of $p$-power order that are finite and flat over $\Spec R$, and where morphisms are homomorphisms of $R$-group schemes.
\end{definition}

\begin{definition}[Tate twisting]  \label{tate-twisting}
  Let $H$ be an object of $\fin^{p}_{R}$ and suppose that $n$ is a power of $p$ such that $nH = 0$. Write $(-)^{*}$ for Cartier duality, and write $(-)^{\wedge}$ for Pontryagin duality on the full subcategory of \'etale objects on $\fin^{p}_{R}$. For $H$ \'etale, define $H(1) := (H^{\wedge})^{*}$. For $H$ multiplicative, define $H(-1) := (H^{*})^{\wedge}$. 
\end{definition}

We will soon make use of the canonical isomorphism $H(1) \simeq H\otimes\mu_{n}$. Recall the morphism
  \[
    \mu: Y\times_{\OK} Y^{D}\longto \Z_{\OK}
  \]
displayed in \eqref{eq:mu-def}.  Let $Y^{D,\vee} := \Hom(Y^{D}, \Z_{\OK})$. This is a twisted constant $\OK$-group. We define a homomorphism
  \begin{align}
    \label{eq:mu1}
    \mu_{Y}: Y\longto Y^{D,\vee}
  \end{align}
  on local sections by $\mu(z):\chi\mapsto \mu(z,\chi)$. Twisting gives a homomorphism of finite $\OK$-groups of multiplicative type
\begin{align}
  \label{eq:twist-of-muY}
Y/nY(1) \longto Y^{D,\vee}/nY^{D,\vee}(1).  
\end{align}
Recall that $Y^{D}$ is defined so that the torus $T$ in $\widetilde A$ has the form $T = \Hom(Y^{D},\Gm)$. Using the isomorphisms of abelian sheaves
\[
  T = \underline{\Hom}(Y^{D},\Z)\otimes \Gm = Y^{D,\vee}\otimes\Gm
\]
and writing $T[n] = \ker(Y ^{D,\vee}\otimes \Gm \overset{n}\to Y^{D,\vee}\otimes \Gm)$, we see that
\begin{align}
  \label{eq:T-n-as-quotient-twist}
    Y^{D,\vee}/nY^{D,\vee}(1)
  = Y^{D,\vee}/nY^{D,\vee}\otimes \mu_{n}
  \simeq Y^{D,\vee}\otimes \mu_{n}
  = T[n].
\end{align}
Composing \eqref{eq:twist-of-muY} and the identification \eqref{eq:T-n-as-quotient-twist} with the canonical closed immersion $T[n]\hooklongrightarrow \widetilde A[n]$ gives a homomorphism.
\begin{align}
  \label{eq:nu-n-def}
\nu_{n}: Y/nY(1) \longto \widetilde A[n].  
\end{align}

\begin{example}
In the case of a Tate curve, where $Y/nY = \Z_{\OK}/n\Z_{\OK}$ and $\widetilde A = \GmOK$, the morphism $\nu_{n}$ is the Cartier dual of the morphism $N_{n}^{*}$ defined at the end of \Cref{subsec:q-tate-eg}.
\end{example}

The following result is needed for the proof of \Cref{main-theorem-torsion}.

\begin{proposition}
  \label{ker-gives-component-group}
  There is a canonical isomorphism
  \[
    \ker \nu_{n}(-1) \simeq \Phi_{A_{K}}[n].
  \]
  \proof{
By \cite{faltings-chai} Chapter III, Corollary 8.2 (see also \cite{grothendieck-monodromy}, 11.5), there is an isomorphism $\Phi_{A_{K}}\simeq \coker \mu_{Y}$, where $\mu_{Y}: Y\to Y^{D, \vee}$ is the morphism of \eqref{eq:mu1}. The groups $Y_{K}(\Kbar)$, $Y^{D, \vee}_{K}(\Kbar)$ are free of rank $\rk T = \rk T^{D}$. Let $c := \#\Phi_{A_{K}}(\Kbar)$.  Then $\mu_{Y}$ gives a homomorphism $c\mu_{Y}: (1/c)Y_{K}(\Kbar)\to Y^{D, \vee}$ which also has cokernel $\Phi_{A_{K}}(\Kbar)$. The map $c\mu_{Y}$ induces a morphism $\mu_{Y,c}(-1)$ fitting in the following commutative diagram of $\Gamma_{K}$-modules having exact rows:
    \[
      \begin{tikzcd}
        & Y_{K}(\Kbar)\ar[d,"\mu_{Y}"] \ar[r]
        & \frac{1}{c} Y_{K}(\Kbar) \ar[d, "c\mu_{Y}"] \ar[r]
        & \frac{1}{c}Y_{K}(\Kbar)/Y_{K}(\Kbar)
        \ar[d, "\mu_{Y, c}(-1)"]
        \\
        &  Y_{K}^{D, \vee}(\Kbar) \ar[r, "c"]
        & Y_{K}^{D, \vee}(\Kbar) \ar[r]
        & Y_{K}^{D, \vee}(\Kbar)/cY_{K}^{D, \vee}(\Kbar)
      \end{tikzcd}.
    \]
Since the homomorphism $\mu_{Y}: Y_{K}(\Kbar)\to Y_{K}^{D, \vee}(\Kbar)$ has finite cokernel $\Phi_{A_{K}}(\Kbar)$, it is injective. Since $(1/c)Y_{K}(\Kbar)$ is free, $c\mu_{Y}$ is also injective. The Snake Lemma gives an exact sequence
    \[
      0 \longto 0 \longto \ker\mu_{Y, c}(-1) \longto \Phi_{A_{K}}(\Kbar)
      \overset{c}\longto \Phi_{A_{K}}(\Kbar)\longto \coker \mu_{Y,c}(-1).
    \]
    Since $c$ is the order of $\Phi_{A_{K}}(\Kbar)$, we have $\ker\mu_{Y, c}(-1) = \Phi_{A_{K}}(\Kbar )$. The proposition follows upon passing to the $n$-torsion subgroups. Specifically, restricting $\mu_{Y,c}(-1)$ to $\frac{1}{n}Y_{K}(\Kbar)/Y_{K}(\Kbar)$, we obtain a morphism whose image is killed by $n$. This induces a morphism $(1/n)Y_{K}/Y_{K} \to Y_{K}^{D, \vee}/nY_{K}^{D, \vee}$. Making the identification $(1/n)Y_{K}/Y_{K} = Y_{K}/nY_{K}$, we see that this is exactly the morphism $\nu_{n}(-1)$.
    \qed}
\end{proposition}

\section{Maximal 1-crystalline submodules}
In this section we develop results from integral $p$-adic Hodge theory, and we define the functor $\Crys_{r}$ of maximal $r$-crystalline submodules and its continuous derived functors. Beyond some innovations in presentation, including our use of the result \Cref{liu-breuil-conjecture}, which was not available to Kim and Marshall when they wrote \cite{kim-marshall}, we have added little to the results of \cite{kim-marshall}, which we recommend (along with Marshall's more expansive thesis \cite{marshall-thesis}) as a reference.

Throughout the section, we adopt the notation of \Cref{subsec:notation}.

\subsection{Categories of representations and full faithfulness theorems}

To define the continuous derived functors $R^{i}\Crys_{r}$ of the functor of maximal $r$-crystalline sub-representations, following \cite{kim-marshall} and \cite{marshall-thesis}, we first define various categories of $\Gamma_{K}$-modules that can be constructed from $\Zp$-lattices in crystalline $\Qp$-representations. We do this now, but we  do not review Fontaine's theory of period rings, which are described in many sources (for example, \cite{brinon-conrad}). We normalize Hodge-Tate weights so that the $p$-adic cyclotomic character has 1 as its only Hodge-Tate weight.

\begin{definition}
  \label{def:cats-of-ast-modules}
  Let $r\in \Z_{\geq 0}$. Let $\ast\in\{\dR, \stbl, \crys\}$.
  \begin{enumerate}[(i)]
  \item Let $\ModZpGK$ denote the category of topological $\Zp[\Gamma_{K}]$-modules, i.e., the category of topological spaces $M$ equipped with a continuous morphism $\Zp[\Gamma_{K}]\times M\to M$ giving $M$ the structure of a $\Zp[\Gamma_{K}]$-module.
  \item Let $\RepZp$ denote the full subcategory of $\ModZpGK$ consisting of the objects that are finite and free over $\Zp$ and that are equipped with the $p$-adic topology.
  \item Let $\ModtorsZpGK$ denote the full subcategory of $\ModZpGK$ consisting of torsion objects with the discrete topology.
  \item Let $\RepastrQp$ denote the full subcategory of the $B_{\ast}$-admissible representation of $\Gamma_{K}$ on finite-dimensional $\Qp$-vector spaces having all Hodge-Tate weights in the interval $[0,r]$.
  \item Let $\Repastrfree$ denote the full subcategory of $\RepZp$ consisting of those objects $L$ such that $L[1/p] := L\otimes_{\Zp} \Qp$ is an object of $\RepastrQp$.
  \item Let $\Repastrtors$ denote the full subcategory of $\ModtorsZpGK$ consisting of those objects $V$ such that there exists an object $L$ of $\Repastrfree$ and a surjective morphism $L \to  V$ in the category $\ModZpGK$.
  \item Let $\Repastrft$ denote the full subcategory of $\ModZpGK$ consisting of those $M$ such that $M$ is finite type over $\Zp$ and $M$ is isomorphic to the inverse limit of a system of objects of $\Repastrtors$.
  \end{enumerate}
\end{definition}

\begin{example}
  Any unramified representation is crystalline with $0$ as its only Hodge-Tate weight. In particular, the trivial character $\Gamma_{K}\to \{\id\}\subset \GL_{1}(\Qp)$ is crystalline. Conversely, any de Rham representation with $0$ as its only Hodge-Tate weight is unramified (see \cite{brinon-conrad} 8.3.5).
\end{example}

\begin{example}
  The $p$-adic cyclotomic character $\Zp(1)$ is 1-crystalline, hence each of its torsion quotients, $\Z/p^{m}\Z(1)$ is 1-crystalline. In general, any crystalline $\Qp$-character of $\Gamma_{K}$ is the tensor product of an unramified $\Qp$-character with an integer power of $\chi$. In particular, a character of $\Gamma_{K}$ is crystalline if and only if it is de Rham.
\end{example}

\begin{example}
  Let $A_{K}$ be an abelian variety over $K$. Then \[V_{p}A_{K} := T_{p}A_{K}\otimes_{\Zp}\Qp\]
  is a de Rham representation of dimension $2\dim (A_{K})$ having Hodge-Tate weights in $[0,1]$. When $A_{K}$ has semistable reduction, $V_{p}A_{K}$ is semistable. This fact, which is well known, follows from \Cref{BCC-main-theorem} and \Cref{monodromy-pushout-is-1sst}. When $A_{K}$ has good reduction, $V_{p}A_{K}$  is crystalline; see \Cref{coleman-iovita-breuil}.
\end{example}

The following classification result summarizes work of many people completed over the course of fifty years.

\begin{theorem}[Tate; Raynaud; Kisin; Kim; see \Cref{1-crys-and-geom-attributions}]
  \label{1-crys-things-are-gen-fibers}
  One has the following relationships between $1$-crystalline $\Gamma_{K}$-modules and finite flat $\OK$-groups. Let $\pdiv_{\OK}$ denote the category of $p$-divisible groups over $\OK$.
  \begin{enumerate}[(i)]
  \item The Tate module of any object of $\pdiv_{\OK}$ is an object of $\RepcrysoneZp$. 
  \item The generic fiber of any object of $\fin^{p}_{\OK}$ is an object of $\Repcrysonetors$.
  \item Every object of $\RepcrysoneZp$ is the Tate module an object of $\pdiv_{\OK}$.
  \item Every object of $\Repcrysonetors$ is the generic fiber of an object of $\fin^{p}_{\OK}$.
  \end{enumerate}
\end{theorem}

\begin{remark}
  \label{1-crys-and-geom-attributions}
  Part (i) was proved by Tate in different terms in \cite{tate-p-divisible}, Section 4, Corollary 2. Part (ii) is attributed to Oort in \cite{brinon-conrad}, but it also follows from Raynaud's stronger result that every object of $\fin^{p}_{\OK}$ is the kernel of an isogeny of abelian $\OK$-schemes (see \cite{berthelot-breen-messing-II}, 3.1.1). Parts (iii) and (iv) were proven by Kisin in \cite{kisin-f-crystals} for $p>2$ and for $p=2$ in some cases (see also \cite{kisin-shimura-models}), and was proved unconditionally by Wansu Kim in the case for $p=2$. See \cite{kim-2-adic} 2.3 for an explanation of Kisin's work, and 4.1 for the extension to the case $p=2$.
\end{remark}

We are primarily concerned with objects of the categories $\Repastrfree$ and $\Repastrtors$. A key relationship between these was conjectured by Breuil and proved by Liu.

\begin{theorem}[\cite{liu-torsion-fontaine}, 7.2 and 7.3]
  \label{liu-breuil-conjecture}
  Let $r\in \Z_{\geq 1}$. Let $\ast\in \{\stbl, \crys\}$. Let $L$ be an object of $\RepZp$. Then $L$ is an object of $\Repcrysrfree$ (resp., $\Repsstrfree$) if and only if, for each $m\geq 1$, $L/p^{m}L$ is an object of $\Repcrysrtors$ (resp.~ $\Repsstrtors$).
\end{theorem}

Let $L\in \Repsstrfree$ and assume $L$ is not crystalline. By \Cref{liu-breuil-conjecture}, we know that $L/p^{m}L$ is not an object of $\Repcrysrtors$ for all $m$ sufficiently large. However, it can happen that $L/p^{m}L$ is $r$-crystalline for small values of $m$, as we see in \Cref{tate-curve-prop}.

Recall from \Cref{subsec:notation} the subgroup $\Gamma_{\varpi^{\flat}}$ of $\Gamma_{K}$. Let $\Rep_{\free}(\Gamma_{\varpi^{\flat}})$ (resp., $\Rep_{\tors}(\Gamma_{\varpi^{\flat}})$) denote the category of finite-type free (resp., torsion) $\Zp$-modules equipped with an action of $\Gamma_{\varpi^{\flat}}$. The following result is essential in the study of crystalline $\Qp$-representations.

\begin{proposition}[\cite{kisin-f-crystals} \cite{beilinson-tavares}; \cite{breuil-integral} 3.4.3 in the case $r=1$]
  For any $r\in \Z_{\geq 0}$, the restriction functor
$      \Res^{\Gamma_{K}}_{\Gamma_{\varpi^{\flat}}}:
      \Repcrysrfree
      \longto
      \Rep_{\free}(\Gamma_{\varpi^{\flat}})
$ 
  is fully faithful.
\end{proposition}

In general, the analogous statement for torsion $\Gamma_{K}$-modules fails for reasons which are related to the bound in Raynaud's work on prolongations.

\begin{proposition}[\cite{ozeki-equivariant} 1.2; others in special cases, e.g. $r=1$, $p>2$ case by \cite{breuil-integral} 3.4.3. See the introduction to \cite{ozeki-equivariant} for a detailed history of this result.]
  \label{torsion-ff}
  Let $r \in \Z_{\geq 0}$ and let $e$ be the ramification index of $K/\Qp$. Suppose $e(r-1)<p-1$. Then the restriction functor
  \[
\Res^{\Gamma_{K}}_{\Gamma_{\varpi^{\flat}}}: \Repcrysrtors \longto \Reptors(\Gamma_{\varpi^{\flat}})
  \]
  is fully faithful.
\end{proposition}

The utility of this result was illustrated in \Cref{subsec:q-tate-eg}.

In \Cref{ast-standard-results} and \Cref{ast-ftype-compats} we record some basic closure properties of the categories defined in \Cref{def:cats-of-ast-modules}. We omit the straightforward proofs, most of which have already appeared in detail in \cite{marshall-thesis}, Section 2.2.2.

\begin{proposition}
  \label{ast-standard-results}
  Let $\sC$ be any of the categories defined in  \Cref{def:cats-of-ast-modules}. By definition, $\sC$ is a full subcategory of $\ModZpGK$.
\begin{enumerate}[(a)]
\item Any $\ModZpGK$-subobject of an object of $\sC$ is also an object of $\sC$. 
  \item Suppose $\sC$ is one of the categories containing a nontrivial torsion object.  If $M$ is a $\ModZpGK$-quotient of an object of $\sC$, then $M$ is an object of $\sC$.
  \item If $M$ is a $\ModZpGK$-direct sum of finitely-many objects of $\sC$, then $M$ is an object of $\sC$.
\end{enumerate}
\end{proposition}

\begin{lemma}
  \label{ast-ftype-compats}
  Let $\ast\in \{\stbl, \crys\}$. Let $M$ be an object of $\Repastrft$ and let $M'$ be a $\ModZpGK$-subobject of $M$.
  \begin{enumerate}[(i)]
  \item If $M'$ is free, then $M'$ is an object of $\Repastrfree$.
  \item If $M'$ is torsion, then $M'$ is an object of $\Repastrtors$.
  \item Write $M = M_{\free}\oplus M_{\tors}$ for the decomposition of $M$ as a direct sum of a free $\Zp$-module and its maximal torsion $\Zp$-submodule. Then $M_{\free}$ (resp., $M_{\tors}$) is an object of $\Repastrfree$ (resp., $\Repastrtors$).
   \end{enumerate}
\end{lemma}

Our proof of the key result \Cref{monodromy-criterion} relies on \Cref{baer-ht-closure} below. We make one preliminary remark. Recall that the set $\Ext(B,A)$ of isomorphism classes of Yoneda extensions in an Abelian category $\sA$ becomes a group under the Baer sum. Given objects $E_{1}$ and $E_{2}$ of $\sA$ realizing two objects of $\Ext(B,A)$, we remark that the sum object can be explicitly constructed by taking a subquotient of $E_{1}\times E_{2}$  (see \StacksTag{010I} and  \StacksTag{05PJ}).

\begin{proposition}
  \label{baer-ht-closure}
  Let $\ast\in \{\dR, \stbl, \crys\}$. Suppose we have two short exact sequences
  \begin{align*}
 \tikzses{\eta(V_{i})}{V'}{\alpha_{i}}{V_{i}}{\beta_{i}}{V''} \qquad i=1,2
  \end{align*}
of objects of $\ModtorsZpGK$. Let $V_{+}$ be the middle term in the Baer sum $\eta(V_{1}) + \eta(V_{2})$. The following statements hold.
\begin{enumerate}[(i)]
\item If $V_{1}$ and $V_{2}$ are objects of $\Repastrtors$, then so is $V_{+}$.
\item If any two of $V_{1}$, $V_{2}$ and $V_{+}$ is an object of $\Repcrysrtors$, then every one of these is an object of $\Repcrysrtors$.
\end{enumerate}
\proof{  
By \Cref{ast-standard-results}, $\Repastrtors$ is stable under finite direct sums and passage to subquotients. Therefore, since $U_{1}$ and $U_{2}$ are visibly $\Gamma_{K}$-stable inside of $V_{1}\oplus V_{2}$, they are both objects of $\Repastrtors$, so $V_{+}$ is in $\Repastrtors$.

Taking a Baer difference, (ii) follows from (i).
\qed }
\end{proposition}

\subsection{Maximal crystalline subobjects and continuous derived functors}
\label{sec:crys-r-derived-functors}

Following Marshall's thesis \cite{marshall-thesis}, we define functors $\Crys_{r}$ and construct their derived functors. As we shall see, the derived functors are constructed using inverse systems of torsion $\Gamma_{K}$-modules. As was noted in \cite{marshall-thesis}, this technique for defining derived functors has its origin in Jannsen's work \cite{jannsen-continuous-etale}. We first define maximal $r$-crystalline subobjects and a category where such things lie. \Cref{sec:crys-r-derived-functors} contains nothing original.

\begin{definition}
  \label{def:crys-r}
  For $M$ in $\ModZpGK$ and for $r\in \Z_{\geq 0}$, define $\Crys_{r}(M)$ to be the sum in $M$ of all $C\subset M$ such that $C$ is an object $\Repcrysrft$.

  Let $\ModcrysrZpGK$ denote the full subcategory of $\ModZpGK$ consisting of those $M$ such that $M$ is isomorphic to a filtered colimit of objects of $\Repcrysrft$.
\end{definition}

\begin{proposition}
  \label{crys-r-props}
  Let $M$ be an object of $\ModZpGK$. Let $r\in \Z_{\geq 0}$.
  \begin{enumerate}[(i)]
  \item The object $\Crys_{r}(M)$ of $\ModZpGK$ is an object of $\ModcrysrZpGK$.
    \item If $M$ is finite and is given the discrete topology, then $\Crys_{r}(M)$ is an object of $\Repcrysrtors$.
  \item If $M$ is an object of $\RepZp$, then $\Crys_{r}(M)$ is an object of $\Repcrysrfree$.
    \item If $M$ is an object of $\RepQp$, then $\Crys_{r}(M)$ is an object of $\RepcrysrQp$.
  \end{enumerate}
  \proof{
Statement (i) is verified by noting that internal sums are filtered colimits. When $M$ is finite type, the $\Zp$-submodule $\Crys_{r}(M)\subset M$ is finite type and torsion with the discrete topology with $\Zp$-rank bounded by that of $M$. It is therefore equal to one of the summands defining $\Crys_{r}(M)$ as a sum of finite type submodules of $M$. If $M$ is torsion (resp., free), then $\Crys_{r}(M)$ is torsion (resp., free). Therefore, (ii) and (iii) follow from \Cref{ast-ftype-compats}, and (iv) follows from (iii).
    \qed}
\end{proposition}

\begin{example}
For a Tate curve $E_{q,K}$, $\Crys_{1}(E_{q,K}[n])$ was described in \Cref{tate-curve-prop}.
\end{example}

\begin{proposition}
  \label{crys-r-functor}
  Let $r\in \Z_{\geq 0}$. There exists a functor \[\Crys_{r}: \ModZpGK \to \ModcrysrZpGK\] that sends an object $M$ of $\ModZpGK$ to $\Crys_{r}(M)$ as defined in \Cref{def:crys-r}. This functor is right adjoint to the forgetful functor $\ModcrysrZpGK\to \ModZpGK$, and hence it is left exact and commutes with limits.
  \proof{
    We have shown in \Cref{crys-r-props} that $\Crys_{r}$ is well-defined on objects. For any morphism $f: M \to M'$ of objects of $\ModZpGK$, define $\Crys_{r}(f) := f|_{\Crys_{r}(M)}. $ It is immediate from the definition of $\Crys_{r}(M')$ that $f|_{\Crys_{r}(M)}$ factors through $\Crys_{r}(M')$.
\qed}
\end{proposition}

\begin{remark}
  Let $L\in \RepZp$. Then $\Crys_{r}(L[1/p])\cap L$ is a $\Zp$-lattice in $\Crys_{r}(L[1/p])$ and $\Crys_{r}(L)[1/p] = \Crys_{r}(L[1/p])$.
\end{remark}

We wish to study the derived functors of $\Crys_{r}$, but to define these using the usual formalism of Grothendieck we must restrict $\Crys_{r}$ to some category with enough injectives. Since injective objects must be $p$-divisible, to compute the derived functors of $\Crys_{r}$ on $L\in \RepZp$, we are lead to instead consider $L\otimes \Qp/\Zp$.

\begin{proposition}
  \label{crys-and-colim}
  Let $L\in \RepZp$. Then $L\otimes \Qp/\Zp$ is an object of $\ModtorsZpGK$. The natural morphism
  \[
    \colim_{m\geq 1} \Crys_{r}(L/p^{m}L)\longto
    \Crys_{r}(L\otimes\Qp/\Zp),
  \] where the colimit is with respect to maps $L/p^{m}L\to L/p^{m+1}L$ given by multiplication by $p$, is an isomorphism.
  \proof{By definition of $\Crys_{r}$, we have $\Crys_{r}(L\otimes\Qp/\Zp) = \colim M_{i}$ with $M_{i}\in \Repcrysrft$. Since $L\otimes \Qp/\Zp$ is torsion and each $M_{i}$ is of finite type and $r$-crystalline, we see that for each $i\in I$ there is some $k(i)\in \Z_{\geq 1}$ such that $M_{i}\to \Crys_{r}(L\otimes \Qp/\Zp)$ factors through $\Crys_{r}(L/p^{k(i)}L)$. This induces an inverse to the map $\colim_{m\geq 1} \Crys_{r}(L/p^{m}L) \longto \Crys_{r}(L\otimes\Qp/\Zp)$.
\qed}
\end{proposition}

Following Kim-Marshall and Jannsen, we work with inverse systems of torsion $\Gamma_{K}$-modules, since the category of such inverse systems has enough injectives.

\begin{proposition}
  \label{cats-have-enough-injectives}
  The categories $\ModZpGK$, $\ModtorsZpGK$, and the category of $\Inv\ModtorsZpGK$ of inverse systems of objects of $\ModtorsZpGK$ each have enough injectives.
  \proof{
The claim for $\ModZpGK$ follows from \StacksTag{04JE}. Let $M$ be an object of $\ModtorsZpGK$. Let $M\to I$ be an injection of $M$ into an injective object of $\ModZpGK$. Since $M$ is torsion, this injection factors through $I_{\tors}$, so we need only show that $I_{\tors}$ is an injective object in the category $\ModtorsZpGK$. Suppose we have morphisms $f:T\to I_{\tors}$ and $g:T\to T'$ in $\ModtorsZpGK$. Then $f$ and $g$ are also morphisms in the category $\ModZpGK$, so we get a morphism $f':T'\to I$ such that $f = f'\circ g : T\to I$. But $T'$ is torsion, so $f'$ factors through $I_{\tors}$. This shows that $I_{\tors}$ is injective, finishing the proof of the claim for $\ModtorsZpGK$. The claim for $\Inv\ModtorsZpGK$ is immediate from what we have.
    \qed}
\end{proposition}

\begin{definition}
  Let $M$ be an object of $\ModtorsZpGK$. Define objects of $\Inv\ModtorsZpGK$ as follows.
  \begin{enumerate}[(i)]
  \item $T(M) := (T_{m})_{m\geq 1}$ with $T_{m} = M[p^{m}]$ and transition maps $T_{m}\to T_{m'}$, for $m'\leq m$, given by multiplication by $p^{m-m'}$.
  \item $V(M) := (V_{m})_{m\geq 1}$ with $V_{m} = M$ and transition maps $V_{m}\to V_{m'}$, for $m'\leq m$, given by multiplication by $p^{m-m'}$.
  \item $\underline{M} := (M_{m})_{m\geq 1}$ the constant system with $M_{m} = M$ and each transition morphism is the identity.
  \end{enumerate}
\end{definition}

\begin{lemma}
  \label{crys-of-Qp-vs-tors-free}
Let $L\in \RepZp$, $j\geq 0$, $r\geq 0$. Then $R^{i}\Crys_{r}(L[1/p])$ has no $p$-torsion.
  \proof{
 Every element of $V(L \otimes \Qp/\Zp)$ is sequence of the form $(a_{j}/p^{m+j})_{j\geq 0}$, where $m\in \Z_{\geq 1}$, $a_{j}\in L$ for all $j\in \Z_{\geq 0}$, and $a_{j+1} \equiv a_{j}$ modulo $p^{m+j}$. From this we see that the morphism $p: V(L \otimes \Qp/\Zp) \to V(L \otimes \Qp/\Zp)$ is an isomorphism. Taking the long exact sequence associated to this isomorphism using the usual derived functor cohomology, we find that $p$ is an isomorphism on $R^{i}\Crys_{r}(V(L \otimes \Qp/\Zp))$ for each $j\geq 0$. In particular, $R^{i}\Crys_{r}(L[1/p])$ has no $p$-torsion.
  \qed}
\end{lemma}

\begin{theorem}[\cite{marshall-thesis}, Proposition 3.2]
  \label{les-exists}
  Let $L$ be an object of $ \RepZp$ and let $r\in \Z_{\geq 0}$. Then there are canonical isomorphisms
\begin{align*}
  \lim\Crys_{r}(T(L \otimes\Qp/\Zp))
  &\simeq \Crys_{r}(L)
  \\
  \lim\Crys_{r}(V(L \otimes\Qp/\Zp))
  &\simeq \Crys_{r}(L[1/p])
  \\
    \lim\Crys_{r}\left(\underline{L \otimes\Qp/\Zp}\right)
  &\simeq \Crys_{r}(L\otimes\Qp/\Zp).
\end{align*}
For $i\geq 1$, define
\begin{align*}
  R^{i}\Crys_{r}(L)
  &:= R^{i}(\lim\Crys_{r})(T(L \otimes\Qp/\Zp))
  \\
  R^{i}\Crys_{r}(L[1/p])
  &:= R^{i}(\lim\Crys_{r})(V(L \otimes\Qp/\Zp))
  \\
  R^{i}\Crys_{r}(L\otimes \Qp/\Zp)
  &:= R^{i}(\lim\Crys_{r})\left(\underline{L \otimes\Qp/\Zp}\right).
\end{align*}
The morphisms
\[
  \pi_{m} = \id\otimes p^{m}: L \otimes \Qp/\Zp \longto L\otimes \Qp/\Zp
\]
defined for $m\in \Z_{\geq 1}$ give rise to a short exact sequence of objects of $\Inv\ModtorsZpGK$
  \begin{align}
    \label{eq:ses-of-inv-systems}
    \ses{T(L\otimes \Qp/\Zp)}{}{V(L\otimes \Qp/\Zp)}{}{\underline{L\otimes \Qp/\Zp}},
  \end{align}
  hence also a long exact sequence in $\ModcrysrZpGK$
  \begin{equation*}
  \begin{tikzcd}[column sep=small]
    0 \ar[r]
  & \Crys_{r}(L) \ar[r]
  & \Crys_{r}(L[1/p]) \ar[r]
  & \Crys_{r}(L\otimes\Qp/\Zp) \ar[lld] \\
  & R^{1}\Crys_{r}(L) \ar[r]
  & R^{1}\Crys_{r}(L[1/p]) \ar[r]
  & R^{1}\Crys_{r}(L\otimes\Qp/\Zp) \ar[lld] \\
  & R^{2}\Crys_{r}(L) \ar[r]
  & \cdots
  & 
\end{tikzcd}
\end{equation*}
\proof{The first three claimed equalities are immediately verified, where for the second we use the fact from \Cref{crys-r-functor} that $\Crys_{r}$ preserves limits.
  
We now show that the morphisms $\pi_{m}$ give rise to a short exact sequence \eqref{eq:ses-of-inv-systems}. Let $i = (i_{m})_{m\geq 1} :T(L\otimes \Qp/\Zp) \to V(L\otimes \Qp/\Zp)$ be the morphism of inverse systems obtained by taking $i_{m}: p^{-m}L/L = (L\otimes \Qp/\Zp)[p^{m}] \to L\otimes \Qp/\Zp$ to be the inclusion. Then the cokernel of $i_{m}$ is isomorphic $L\otimes \Qp/\Zp$, so we have a short exact sequence
  \[\ses{\frac{1}{p^{m}}L/L}{[p^{m}]\circ i_{m}}{L \otimes \Qp/\Zp}{\pi_{m}}{L \otimes \Qp/\Zp}\]
  The transition maps in the systems $T(L\otimes \Qp/\Zp)$ and $V(L\otimes \Qp/\Zp)$ are both multiplication by $p$, so they commute with the maps $[p^{m}]\circ i_{m}$. A direct verification shows that $\pi_{m+1} = p\circ\pi_{m}$. Therefore, we have the short exact sequence \eqref{eq:ses-of-inv-systems}. As we saw in \Cref{cats-have-enough-injectives}, the category $\Inv\ModtorsZpGK$ has enough injectives. With this, the existence of the long exact sequence is standard.
\qed}
\end{theorem}

\section{Finite flat group schemes, monodromy, and 1-semistable modules}
\label{sec:ffgs-monodromy-1sstmodules}

In this section, we develop the theory of finite flat group schemes with monodromy structures, which was initiated by initiated by Kato (see, e.g., \cite{kato-log-dieudonne}). For the relationship between this perspective and the perspective of log finite flat group schemes, the interested reader should refer to the introduction \Cref{sec:introduction} and Remarks \ref{logarithmic-enhancements}, \ref{log-ffgss} and \ref{kato-tate-torsion}. After discussing some standard results and describing some important examples, we introduce a category $\ExtMonR$ that provides a useful, flexible generalization of the category $\fin^{p,N}_{\OK}$ introduced by Kato. We then develop results on Kato's generic fiber functor. This part of the work amounts to making explicit calculations with Baer sums and certain pushouts. We then describe the relationship to the theory of degenerations of abelian varieties and 1-motives developed in \Cref{sec:monodromy}. In the final part of the section, we describe a monodromy criterion for 1-crystallinity. This result uses the Torsion Full Faithfulness Theorem, which can also be viewed as an explicit calculation involving certain Kummer 1-cocycles of $\Gamma_{K}$

Throughout this section, we adopt the notation of \Cref{subsec:notation}.

\subsection{Finite flat group schemes with monodromy}
\label{subsec:fin-and-finN}

Until \Cref{def:nu-pushout}, let $R$ be a field or a Henselian discrete valuation ring. It is well known that $\fin^{p}_{R}$ (\Cref{def:finpR}) inherits the structure of an exact category by embedding it in $\ShAbfppf{R}$. With this structure, a short exact sequence of objects of $\fin^{p}_{R}$ is a sequence
\begin{align*}
      \ses{\eta^{-1}}{i}{\eta^{0}}{\pi}{\eta^{1}}
  \end{align*}
  in $\ShAbfppf{R}$, where $i$ is given by a closed immersion of $R$-schemes and $\pi$ is given by a faithfully flat morphism of $R$-schemes.
  
The full subcategory of $\fin^{p}_{R}$ consisting of the \'etale objects is equivalent to the category of continuous representations of the \'etale fundamental group of $\Spec R$ on discrete torsion $\Zp$-modules. In particular, since $K$ has characteristic 0, $\ModtorsZpGK$ is equivalent to $\fin^{p}_{K}$. However, it is known (for instance, by \Cref{liu-breuil-conjecture}) that not every object of $\fin^{p}_{K}$ is the generic fiber of an object of $\fin^{p}_{\OK}$. The objects of  $\fin^{p}_{K}$ that do prolong over $\OK$ have traditionally been called ``finite flat''. As we saw in \Cref{1-crys-things-are-gen-fibers}, the condition of being finite flat in this sense is the same the condition of being 1-crystalline.
  
  \begin{example}
  Suppose $G$ is a semiabelian $R$-scheme that sits in a short exact sequence
  \begin{align}
    \label{eq:TGB-ses}
    \ses{T}{}{G}{}{B}
  \end{align}
 where $T$ is an $\OK$-torus and $B$ is an abelian $\OK$-scheme. It can be shown (\cite{grothendieck-montreal} Example III 6.4) that, for any $n\in \bZ_{\geq 2}$ a power of $p$, there is a short exact sequence in $\fin^{p}_{\OK}$
  \begin{align}
    \label{eq:global-semiab-torsion-ses}
    \ses{T[n]}{}{G[n]}{}{B[n]}.
  \end{align}

  Note that, for a general semiabelian $R$-scheme, the torsion $G[n]$ need not be an object of $\fin^{p}_{R}$. For $G$ as in \eqref{eq:TGB-ses}, however, the finite flat groups $G[n]$ fit into a $p$-divisible group $G[p^{\oo}]$.
  
  In the case where $G$ is the Raynaud extension $G = \widetilde A$ attached to a $K$-variety with semistable reduction, we have the Tate module $T_{p}\widetilde A_{K}$. This will play an important role in \Cref{sec:main-results}.
\end{example}

\begin{example}
  \label{R-1-motive-ses}
  For $\sM$ an $R$-1-motive and $n\in \Z_{\geq 2}$ a power of $p$, we have the short exact sequence $\eta(\sM, n)$ in $\fin^{p}_{R}$ as shown in \eqref{eq:eta-M-n}.
\end{example}

In addition to the short exact sequences of the form \eqref{eq:eta-M-n}, we will also work with the connected-\'etale sequences.

\begin{definition}[connected-\'etale sequences 
  ]
  \label{ceseq-def}
    For $H$ an object of $\fin^{p}_{R}$, denote by $\ceseq(H)$ the functorial connected-\'etale sequence
    \[
      \tikzses{\ceseq(H):}{H^{\circ}}{}{H}{}{H^{\et}}.
    \]
\end{definition}

\begin{remark}
  In \Cref{R-1-motive-ses}, though $Y/nY$ is \'etale, the exact sequence $\eta(\sM,n)$ need not be the connected-\'etale sequence of $\sM[n]$. We can say more. Taking \'etale parts is an exact functor. Applying this functor to \eqref{eq:global-semiab-torsion-ses} and noting that $T[n]$ is connected, we find that $G[n]^{\et} = B[n]^{\et}$. Since $Y/nY$ is \'etale, taking \'etale parts in $\eta(\sM, n)$ gives a canonical short exact sequence
\begin{align}
  \label{eq:R-1-motive-et-ses}
  \ses{B[n]^{\et}}{}{\sM[n]^{\et}}{}{Y/nY}.
\end{align}
It follows that $\eta(\sM, n)$ is the connected-\'etale sequence of $\sM[n]$ if and only if $B[n]^{\et}$ is trivial.
\end{remark}

We now enrich our finite flat group schemes with monodromy structure.

\begin{definition}
  Let $R$ be as in \Cref{def:finpR}.
  Let $\ExtMonR$ denote the category of pairs $(\eta^{\bullet}, \nu)$ such that
  \begin{enumerate}[(i)]
  \item $\eta^{\bullet} = [\eta^{-1}\to \eta^{0}\to \eta^{1}]$ is a short exact sequence of objects of $\fin^{p}_{R}$ regarded as a cochain complex concentrated in degrees $-1$, $0$ and $1;$
  \item the term $\eta^{1}$ is an \'etale $R$-group; 
  \item $\nu \in \Hom_{\fin^{p}_{R}}(\eta^{1}(1), \eta^{-1})$.
  \end{enumerate}
A morphism $f: (\eta_{1}^{\bullet}, \nu_{1}) \to (\eta_{2}^{\bullet}, \nu_{2})$ of objects of $\ExtMonR$ is a morphism $f: \eta_{1}^{\bullet} \to \eta_{2}^{\bullet}$ such that the following diagram\footnote{Note that $f^{-1}$ is $f^{\bullet}$ with $\bullet=-1$, not the inverse of $f$.} 
  \[
    \begin{tikzcd}
      \eta_{1}^{1}(1)\ar[r, "\nu_{1}"]\ar[d, "f^{1}(1)"]
      & \eta_{1}^{-1} \ar[d, "f^{-1}"] \\
      \eta_{2}^{1}(1)\ar[r, "\nu_{2}"] & \eta_{2}^{-1}
    \end{tikzcd}
  \]
commutes. Recalling our notation for connected-\'etale sequences in \Cref{ceseq-def}, define $\fin^{p, N}_{R}$ to be the category of pairs $(H, N)$ such that $(\ceseq(H), N)$ is an object of $\ExtMonR$, and where a morphism $(H_{1}, N_{1})\to (H_{2}, N_{2})$ of objects of $\fin^{p, N}_{R}$ is the same as a morphism $(\ceseq(H_{1}), N_{1})\to (\ceseq(H_{2}), N_{2})$.
\end{definition}

\begin{remark}
  \label{log-ffgss}
  It is shown in \cite{kato-log-dieudonne} that a choice of uniformizer of $\OK$ identifies $\fin^{p, N}_{\OK}$ with a category of log finite flat $\OK$-group schemes. This category consists of certain representable sheaves on the big Kummer log flat site of the fs log scheme $S^{\log}$ defined in \Cref{logarithmic-enhancements}.
\end{remark}

\begin{remark}
The categories $\ExtMonR$ and $\fin^{p, N}_{R}$ each inherit the structure of an exact category from that of $\fin^{p}_{R}$. There is a fully faithful functor $\fin^{p}_{R}\to \fin^{p, N}_{R}$ given on objects by $H\mapsto (H, 0)$. In the context of \Cref{log-ffgss}, this functor amounts to equipping $H$ with the induced log structure.
\end{remark}

For the rest of \Cref{subsec:fin-and-finN}, we specialize to the case were $R$ is $K$ or $\OK$.

\begin{construction}
  \label{def:nu-pushout}
  Let $(\eta^{\bullet},\nu)\in \ExtMon$, so that $\eta^{\bullet}$ gives a class in $\Ext_{\fin^{p}_{K}}(\eta^{1}_{K}, \eta^{-1}_{K})$. We will construct another class in $\Ext_{\fin^{p}_{K}}(\eta^{1}_{K}, \eta^{-1}_{K})$.

  Let $n$ be the order of $\eta^{1}$, so that $\eta^{1}$ is a $\Z/n\Z_{\OK}$-module. Recall that for any Tate curve $E_{q,K}$, the finite $K$-group of $n$-torsion points $E_{q,K}[n]$ lies in a short exact sequence
\begin{align*}
\tikzses{\theta^{q}_{n,K}:}{\mu_{n,K}}{}{E_{q,K}[n]}{}{\Z/n\Z_{K}}.
\end{align*}
Consider the special case where $q$ is a uniformizer $\varpi$. Tensor the sequence $\theta^{\varpi}_{n,K}$ of free $\Z/n\Z$-modules with $\eta^{1}_{K}$ and push out by $\nu$ to obtain
\begin{align*}
  \MP(\eta^{\bullet}, \nu)_{K} :=
  \nu_{*}(\theta^{\varpi}_{n,K}\otimes \eta^{-1}_{K})
  \in \Ext_{\fin^{p}_{K}}(\eta^{1}_{K}, \eta^{-1}_{K}).
\end{align*}
By definition, $\MP(\eta^{\bullet}, \nu)_{K}$ is the bottom row in the pushout diagram
\[
  \begin{tikzcd}
    0\ar[r]
    & \eta^{1}_{K}(1)\ar[r]\ar[d, "\nu"]
    & \eta^{1}_{K}\otimes E_{\varpi, K}[n]\ar[r]\ar[d]
    & \eta^{1}_{K}\ar[r]\ar[d, equals]
    & 0
    \\
    0\ar[r]
    & \eta^{-1}_{K}\ar[r]
    & \MP(\eta^{\bullet}, \nu)_{K}^{0} \ar[r]
    & \eta^{1}_{K} \ar[r]
    & 0,
  \end{tikzcd}
\]
and
\begin{align}
  \label{eq:nu-star-eta-0-explicit}
  \MP(\eta^{\bullet}, \nu)_{K}^{0}
  = \frac{(\eta^{1}_{K}\otimes E_{\varpi, K}[n])\oplus \eta^{-1}_{K}}{\{(y,0)-(0, \nu(y)):\: y\in \eta^{1}_{K}(1)\}}.
\end{align}
\end{construction}

\begin{remark}
  \label{kato-tate-torsion}
Let $n\in \Z_{\geq 2}$ be a power of $p$ and let $S^{\log}$ be as in \Cref{log-ffgss}. Kato has constructed a log finite flat $S^{\log}$-group that has generic fiber isomorphic to $E_{\varpi, K}[n]$. Using this, one can mimic \Cref{def:nu-pushout} integrally, i.e., in a certain category of fs log schemes over $S^{\log}$. With this in mind, we sometimes carry a subscript of $K$ in our notation even in situations where we do not define an $\OK$-model.
\end{remark}

We will use the following lemma, which holds in much greater generality, in our definition of the generic fiber functor on $\ExtMonR$.

\begin{lemma}
    \label{coch-les-lemma}
    For $i=1,2,3$, let $V_{i}^{\bullet}$ be acyclic cochain complexes of objects of $\fin^{p}_{K}$ concentrated in degrees $-1$, $0$ and $1$. Suppose we have two morphisms
    \begin{align}
      \label{eq:eta-chain-complex}
      V_{3}^{\bullet}
      \overset{f_{2}^{\bullet}}\longto V_{2}^{\bullet}
      \overset{f_{1}^{\bullet}}\longto  V_{1}^{\bullet}
    \end{align}
    giving short exact sequences in degrees $1$ and $-1$. 
    \begin{enumerate}[(i)]
    \item Then $f_{2}^{0}$ is injective and $f_{1}^{0}$ is surjective.
    \item If $\ker f_{1}^{0}\subset \im f_{2}^{0}$, then in fact $\ker f_{1}^{0}= \im f_{2}^{0}$. Therefore, by (i), \eqref{eq:eta-chain-complex} gives a short exact sequence in degree $0$.
    \end{enumerate}
    \proof{
The morphisms of \eqref{eq:eta-chain-complex} fit into the following commutative diagram with exact rows and where the outer two columns are short exact: 
\begin{equation*}
  \begin{tikzcd}
    0 \ar[r]
    & V_{3}^{-1}  \ar[r] \ar[d, "f_{2}^{-1}"]
    & V_{3}^{0}  \ar[r] \ar[d, "f_{2}^{0}"]
    & V_{3}^{1} \ar[r] \ar[d, "f_{2}^{1}"]
    & 0 \\
    0 \ar[r]
    & V_{2}^{-1}  \ar[r] \ar[d, "f_{1}^{-1}"]
    & V_{2}^{0}  \ar[r] \ar[d, "f_{1}^{0}"]
    & V_{2}^{1} \ar[r] \ar[d, "f_{1}^{1}"]
    & 0 \\
    0 \ar[r]
    & V_{1}^{-1}  \ar[r]
    & V_{1}^{0}  \ar[r]
    & V_{1}^{1} \ar[r]
    & 0    
  \end{tikzcd}
\end{equation*}
Applying the Snake Lemma to the top two rows and using the hypotheses, we obtain an exact sequence
\begin{align*}
  0\longto \ker f_{2}^{0} \longto 0 \longto V_{1}^{-1} \longto \coker f_{2}^{0} \longto V_{1}^{1}.
\end{align*}
This shows that $\ker f_{2}^{0} = 0$. Similarly, applying the Snake Lemma to the bottom two rows gives an exact sequence
\begin{align}
  \label{eq:snake-lemma-2}
    V_{3}^{-1} \longto \ker f_{1}^{0} \longto V_{3}^{1} \longto 0 \longto \coker f_{1}^{0} \longto 0,
\end{align}
which shows that $\coker f_{1}^{0} = 0$. This has proven (i).

We prove (ii). Note that the morphism $d_{2}:\ker f_{1}^{-1} \longto \ker f_{1}^{0}$ is injective by assumption, so \eqref{eq:snake-lemma-2} gives a short exact sequence
\begin{align*}
    0\longto V_{3}^{-1} \longto \ker f_{1}^{0} \longto V_{3}^{1} \longto 0.
\end{align*}
We therefore have equations relating orders of finite abelian groups
\[
  \#\ker f_{1}^{0}
  = (\#V_{3}^{1})(\#V_{3}^{-1})
  = \# V_{3}^{0}.
\]

By (i), $f_{2}^{0}$ is injective, so $\#\im f_{2}^{0} = \# V_{3}^{0}$. We see that $\ker f_{1}^{0}$ and $\im f_{2}^{0}$ are both of order $\# V_{3}^{0}$. Therefore, if $\ker f_{1}^{0}\subset \im f_{2}^{0}$, then $\ker f_{1}^{0}= \im f_{2}^{0}$.
      \qed}
  \end{lemma}

\begin{proposition}  
  \label{monodromy-pushout-is-exact}
The association $(\eta^{\bullet},\nu)\mapsto \MP(\eta^{\bullet},\nu)_{K}$ gives rise to an exact functor
\begin{align}
  \label{eq:pushout-functor}
  \ExtMon \longto \Ext_{\fin^{p}_{K}}.
\end{align}
  \proof{
    The universal property of the pushout defining $P_{\varpi}(\eta^{\bullet},\nu)^{0}_{K}$
gives the desired functor \eqref{eq:pushout-functor}. We now show this functor is exact. Let
\begin{equation}
  \label{eq:etas-ses}
  \ses{
    (\eta_{1}^{\bullet}, \nu_{1})}{f_{2}^{\bullet}}{
    (\eta_{2}^{\bullet}, \nu_{2})}{f_{3}^{\bullet}}{
    (\eta_{3}^{\bullet}, \nu_{3})
  }
\end{equation}
be a short exact sequence of objects of $\ExtMon$, which means simply that \eqref{eq:etas-ses} is a short exact sequence of acyclic cochain complexes
\begin{align}
  \label{eq:eta-ses-no-push}
  \ses{\eta_{1, K}^{\bullet}}{f_{2}^{\bullet}}{
       \eta_{2, K}^{\bullet}}{f_{3}^{\bullet}}{
       \eta_{3, K}^{\bullet}}  
\end{align}
that is compatible with the monodromy maps. We must show that the sequence
\begin{equation}
  \label{eq:nu-star-etas-ses}
  \MP(\eta_{1, K}^{\bullet}, \nu_{1})_{K}
  \overset{\MP f_{2, K}^{\bullet}}\longto 
  \MP(\eta_{2, K}^{\bullet}, \nu_{2})_{K}
  \overset{\MP f_{3, K}^{\bullet}}\longto
  \MP(\eta_{3, K}^{\bullet}, \nu_{3})_{K}
\end{equation}
is short exact in each degree $\bullet \in \{-1, 0, 1\}$. For degrees $\bullet\in \{-1, 1\}$, \eqref{eq:nu-star-etas-ses} agrees with \eqref{eq:eta-ses-no-push},
which is exact by the assumption. Therefore, \eqref{eq:nu-star-etas-ses} is short exact in degrees $-1$ and $1$. By \Cref{coch-les-lemma}, to prove short exactness in degree $0$, it is enough to show that $\ker \MP f_{3, K}^{0}\subset \im \MP f_{2, K}^{0}$. We do this now.
Let $n$ be the maximum of the orders of $\eta_{1}^{1}$ and $\eta_{2}^{1}$. Recall from \eqref{eq:nu-star-eta-0-explicit} that, for $i \in \{1,2,3\}$, we have
\begin{align*}
  \MP(\eta_{i}^{\bullet}, \nu_{i})_{K}^{0}
  = \frac{(\eta_{i, K}^{1}\otimes E_{\varpi, K}[n])\oplus \eta_{i, K}^{-1}}{\{(y,0)-(0, \nu_{i}(y)):\: y\in \eta_{i, K}^{1}(1)\}}.
\end{align*}
We will represent the class of a pair $(a,b) \in (\eta_{i,K}^{1}\otimes E_{\varpi, K}[n])\oplus \eta_{i,K}^{-1}$ using brackets as in $[a,b]$. Let $[a_{2}, b_{2}]\in \MP(\eta_{2}^{\bullet}, \nu_{2})_{K}^{0}$. Then $[a_{2}, b_{2}]\in \ker \MP f_{3, K}^{0}$ if and only if
\begin{align}
  \label{eq:pushout-in-kernel-desc}
  ((f_{3}^{1}\otimes \id)(a_{2}), f_{3}^{-1}(b_{2})) = (y_{3}, -\nu_{3}(y_{3})) \text{ for some } y_{3}\in \eta_{3, K}^{-1}(1).
\end{align}
Suppose \eqref{eq:pushout-in-kernel-desc} holds. We wish to show that $[a_{2}, b_{2}]\in \im \MP f_{2, K}^{0}$, i.e., we wish to show that there exists a pair $(a_{1}, b_{1})\in (\eta_{1, K}^{1}\otimes E_{\varpi, K}[n]) \oplus \eta_{1, K}^{-1}$ and some $y_{2}\in \eta_{2, K}^{1}(1)$ such that 
\begin{align}
  \label{eq:im-desideratum-pushforward}
  (a_{2}, b_{2})
  = ((f_{2, K}^{1}\otimes \id)(a_{1}), 
  f_{2, K}^{-1}(b_{1})) + (y_{2}, -\nu_{2}(y_{2})).  
\end{align}
Choose a lift $y_{2}$ of $y_{3}$ under the surjection $f_{3, K}^{1}(1): \eta_{2, K}^{1}(1) \to \eta_{3, K}^{1}(1)$. Then, by \eqref{eq:pushout-in-kernel-desc}, we have
\[
  ((f_{3, K}^{1}\otimes \id)\times f_{3, K}^{-1})\big(
  (a_{2}, b_{2}) - (y_{2}, -\nu_{2}(y_{2}))
  \big) = (0,0) \in (\eta_{3, K}^{1}\otimes E_{\varpi, K}[n]) \oplus \eta_{3, K}^{-1}.
\]
By exactness of \eqref{eq:eta-ses-no-push} in degrees $-1$ and $1$ and exactness of tensoring against the free $\Z/n\Z$-module $E_{\varpi, K}[n]$, there exists 
$  (a_{1}, b_{1})\in (\eta_{1, K}^{1}\otimes E_{\varpi, K}[n]) \oplus \eta_{1, K}^{-1}$
such that we have equalities
\begin{align*}
  (f_{2, K}^{1}\otimes \id)(a_{1})
  &=  a_{2}-y_{2} \\
  f_{2, K}^{-1}(b_{1})
  &=  b_{2} + \nu_{2}(y_{2}).
\end{align*}
This has shown that we have an equation of the form \eqref{eq:im-desideratum-pushforward}. Since $[a_{2}, b_{2}]$ is arbitrary, we have shown that $\ker \MP f_{3, K}^{0}\subset \im \MP f_{2, K}^{0}$. Using \Cref{coch-les-lemma}, we conclude that the functor \eqref{eq:pushout-functor} is exact.
\qed}
\end{proposition}

\begin{remark}
  \label{monodromy-pushout-is-1sst}
  If $(\eta^{\bullet},\nu)\in \ExtMon$, then $\MP(\eta^{\bullet},\nu)_{K}^{0}$ is an object of $\Repsstonetors$. Indeed, $\eta_{K}^{1}$ is 0-crystalline since it prolongs to an \'etale $\OK$-group scheme, and $E_{\varpi,K}[n]$ can be shown by a direct calculation with the period ring $B_{\stbl}$ to be 1-semistable, so $\eta_{K}^{1} \otimes E_{\varpi,K}[n]$ is 1-semistable. Since $\eta_{K}^{-1}$ is 1-crystalline, being the generic fiber of $\eta^{-1}$, the direct sum $\eta_{K}^{1} \otimes E_{\varpi,K}[n] \oplus \eta_{K}^{-1}$ is 1-semistable, and therefore its quotient $\MP(\eta^{\bullet},\nu)_{K}^{0}$ is 1-semistable.
\end{remark}

We define the generic fiber functor on $\ExtMon$ and show that it is exact.

\begin{proposition}
  \label{generic-fiber-functor-is-exact}
There exists an exact functor
\begin{align}
  \label{eq:gen-fiber-functor-def}
    (-)_{K}: \ExtMon
    \longto \Ext_{\fin^{p}_{K}}
\end{align}
which is given on objects by
  \[
    (\eta^{\bullet}, \nu)_{K}
    := \eta^{\bullet}_{K} + \MP(\eta^{\bullet}, \nu)_{K}.
  \]
  \proof{Let us first address a consequence for notation. Given  $(\eta^{\bullet}, \nu) \in \ExtMon$, the generic fiber $(\eta^{\bullet}, \nu)_{K}$ is another complex concentrated in degrees $-1$, $0$ and $1$. We suppress the index recording these degrees; for instance, the nonzero terms in $(\eta^{\bullet}, \nu)_{K}$ are $(\eta^{\bullet}, \nu)_{K}^{-1}$, $(\eta^{\bullet}, \nu)_{K}^{0}$ and $(\eta^{\bullet}, \nu)_{K}^{1}$, while the expression $(\eta^{0}, \nu)_{K}$ is undefined.
    
 Let $f_{2}^{\bullet}: (\eta_{1}^{\bullet}, \nu_{1}) \longto (\eta_{2}^{\bullet}, \nu_{2})$ be a morphism of objects of $\ExtMon$. We seek to define a morphism
$(\eta_{1}^{\bullet}, \nu_{1})_{K} \longto (\eta_{2}^{\bullet}, \nu_{2})_{K}. $ For $\bullet\in \{1,-1\},$ define $f_{2, K}^{+, \bullet} := f_{2, K}^{\bullet}$. We lighten the notation by setting, for $i=1,2$, $\lambda_{i,K}^{\bullet} := \MP(\eta_{i}^{\bullet},\nu_{i})_{K}$. By the construction of the Baer sum, the term $(\eta_{i}^{\bullet}, \nu_{i})_{K}^{0} = ((\eta_{i}^{\bullet}, \nu_{i})_{K})^{0}$ in degree zero of the complex $(\eta_{i}^{\bullet}, \nu_{i})_{K}$ is 
\begin{align}
  \label{eq:eta-i-baer-desc}
        (\eta_{i}^{\bullet}, \nu_{i})_{K}^{0}
      = \frac{\eta_{i, K}^{0}\times_{\eta_{i, K}^{1}}\lambda_{i,K}^{0}}{
        \im(\Delta_{i}^{-})} ,
\end{align}
where $\Delta_{i}^{-} = (d_{\eta_{i,K}^{\bullet}}(x), -d_{\lambda_{i,K}^{\bullet}}(x))$, and where the morphisms $d_{(-)}$ are the differentials in each of the two complexes.
By functoriality of the Baer sum, the morphism
\[
  f_{2, K}^{\times,0}
  := (f_{2, K}^{0}\times \MP f_{2}^{0}):
  \eta_{1,K}^{0}\times_{\eta_{1, K}^{1}}\lambda_{1, K}^{0}
  \longto
  \eta_{2, K}^{0} \times_{\eta_{2, K}^{1}}\lambda_{2,K}^{0}
\]
induces a morphism
$  f^{+, 0}_{2, K}:
  (\eta_{1}^{\bullet}, \nu_{1})_{K}^{0}
      \to (\eta_{2}^{\bullet}, \nu_{2})_{K}^{0}
$
that lies in a morphism of complexes
\begin{equation}
  \label{eq:varphi-2-K}
    \begin{tikzcd}
      (\eta_{1}^{\bullet}, \nu_{1})_{K}
      \ar[d, "f_{2,K}^{+, \bullet}"]
    & 0\ar[r]
    & \eta^{-1}_{1,K}\ar[r] \ar[d, "f_{2}^{-1}"]
    & (\eta_{1}^{\bullet}, \nu_{1})_{K}^{0} \ar[r]\ar[d, "f^{+, 0}_{2, K}"]
    & \eta^{1}_{1,K} \ar[r]\ar[d,"f_{2}^{1}"]
    & 0
    \\
    (\eta_{2}^{\bullet}, \nu_{2})_{K}
    & 0\ar[r]
    & \eta^{-1}_{2,K}\ar[r]
    & (\eta_{2}^{\bullet}, \nu_{2})_{K}^{0}\ar[r]
    & \eta_{K}^{1}\ar[r]
    & 0
  \end{tikzcd}
\end{equation}

We now show that the functor $(-)_{K}$ is exact. Let
\begin{equation}
  \label{eq:ses-of-complexes-2}
    \ses{
    (\eta^{\bullet}_{1}, \nu_{1})}{f_{2}^{\bullet}}{
    (\eta^{\bullet}_{2}, \nu_{2})}{f_{3}^{\bullet}}{
    (\eta^{\bullet}_{3}, \nu_{3})
  }
\end{equation}
be a short exact sequence of objects of $\ExtMon$. Recall that this means that \eqref{eq:ses-of-complexes-2} is a short exact sequence of cochain complexes
\begin{align}
  \label{eq:eta-ses-no-log-gen-fiber}
\ses{\eta_{1,K}^{\bullet}}{f_{2}^{\bullet}}{\eta_{2,K}^{\bullet}}{f_{3}^{\bullet}}{\eta_{3,K}^{\bullet}}  
\end{align}
that is compatible with the monodromy homomorphisms. We wish to show that
\begin{equation}
  \label{eq:ses-of-complexes-2-K}
  (\eta^{\bullet}_{1}, \nu_{1})_{K}
  \overset{f_{2, K}^{+, \bullet}}\longto
  (\eta^{\bullet}_{2}, \nu_{2})_{K}
  \overset{f_{3, K}^{+, \bullet}}\longto
    (\eta^{\bullet}_{3}, \nu_{3})_{K}
\end{equation}
is short exact in each degree. The sequence \eqref{eq:ses-of-complexes-2-K} is short exact in degrees $-1$ and $1$ by the assumption that \eqref{eq:ses-of-complexes-2} is short exact. By \Cref{coch-les-lemma}, to show \eqref{eq:ses-of-complexes-2-K} is short exact, it suffices to show that $\ker f_{3, K}^{+, 0}\subset \im f_{2, K}^{+, 0}$.

Expand the notation established above for $i=1,2$ to the case $i=3$.  We will represent the class of a pair $(a,b)\in \eta^{0}_{i, K} \times_{\eta^{1}_{i, K}} \lambda^{0}_{i,K}$ in the quotient $(\eta_{i}^{\bullet}, \nu_{i})_{K}^{0}$ by $[a,b]$. Recall that the morphisms $f_{i, K}^{+, 0}$, for $i\in \{2, 3\}$, are defined using the morphisms
    \[
  f_{i, K}^{\times, 0}
  := \big(f^{0}_{i,K}\times \MP f^{0}_{i}\big):
  \eta^{0}_{i-1, K}\times_{\eta^{1}_{i-1, K}}\lambda^{0}_{i-1, K}
  \longto \eta^{0}_{i,K}\times_{\eta^{1}_{i,K}}\lambda^{0}_{i,K}.
    \]

    Let $[a_{2}, b_{2}]\in (\eta_{2}^{\bullet}, \nu_{2})_{K}^{0}$. By the description of $(\eta_{2}^{\bullet}, \nu_{2})_{K}^{0}$ given in \eqref{eq:eta-i-baer-desc}, we have $[a_{2},b_{2}]\in \ker f_{3, K}^{+, 0}$ if and only if there exists $x_{3}\in \eta_{3, K}^{-1}$ such that
    \[
      f_{3, K}^{\times, 0}(a_{2}, b_{2})
      = \big(f^{0}_{3,K}(a_{2}),
      \MP f^{0}_{3}(b_{2})\big)
      = \Delta_{3}^{-}(x_{3}).
    \]
By assumption, \eqref{eq:eta-ses-no-log-gen-fiber}
is short exact in every degree, so there exists $x_{2}\in \eta^{-1}_{2,K}$ such that $f_{3,K}^{-1}(x_{2}) = x_{3}$. Since the diagram
\begin{equation*}
  \begin{tikzcd}
    \eta^{0}_{1,K}\times_{\eta^{1}_{1,K}}
    \lambda^{0}_{1,K}
    \ar[r, "f_{2, K}^{\times, 0}"]
    & \eta^{0}_{2,K}\times_{\eta^{1}_{2,K}}
    \lambda^{0}_{2,K} \ar[r]
    \ar[r, "f_{3, K}^{\times, 0}"]
    & \eta^{0}_{3,K}\times_{\eta^{1}_{3,K}}
    \lambda^{0}_{3,K}
    \\
    \eta_{1,K}^{-1}
    \ar[u, "\Delta_{1}^{-}"]
    \ar[r, "f_{2,K}^{-1}"]
    & \eta_{2,K}^{-1}
    \ar[u, "\Delta_{2}^{-}"]
    \ar[r, "f_{3,K}^{-1}"]
    & \eta_{3,K}^{-1}
    \ar[u, "\Delta_{3}^{-}"]
  \end{tikzcd}
\end{equation*}
commutes, we have
\begin{align*}
  f_{3, K}^{\times, 0}
  \big((a_{2}, b_{2}) -\Delta_{2}^{-}(x_{2})\big)
  &= f_{3, K}^{\times, 0}((a_{2}, b_{2}))
    - f_{3, K}^{\times, 0}(\Delta_{2}^{-}(x_{2}))
  \\
  &= \Delta_{3}^{-}(x_{3}) -
    \Delta_{3}^{-}(f_{3,K}^{-1}(x_{2}))
    \\
  &= \Delta_{3}^{-}(x_{3}) -
    \Delta_{3}^{-}(x_{3}) = 0.
\end{align*}
This has shown that $\big(a_{2}-d_{\eta_{2,K}^{\bullet}}(x_{2}),
  b_{2}+ d_{\lambda_{2,K}^{\bullet}}(x_{2})\big)
  \in \ker f_{3, K}^{\times, 0}$, which is to say that we have
  \begin{align*}
    a_{2}-d^{0}_{\eta_{2,K}}(x_{2})
    &\in \ker f_{3,K}^{0}
    \\
    b_{2}+ d_{\nu_{2,*}\eta_{2,K}}(x_{2})
    &\in \ker \MP f_{3}^{0}.
  \end{align*}
  By exactness of \eqref{eq:eta-ses-no-log-gen-fiber}, there exists $a_{1}\in \eta_{1,K}^{0}$ such that
  \[
    f_{2,K}^{0}(a_{1})
    = a_{2}-d_{\eta_{2,K}^{\bullet}}(x_{2}).
  \]
  By \Cref{monodromy-pushout-is-exact}, there exists $b_{1}\in \lambda_{1,K}^{0}$ such that
  \[
    \MP f_{2,K}^{0}(b_{1})
    = b_{2}+ d_{\lambda_{2,K}^{\bullet}}(x_{2}).
  \]
  We claim that $(a_{1},b_{1})$ gives a pair in the fiber product $\eta^{0}_{1,K} \times_{\eta^{1}_{1,K}} \lambda^{0}_{1,K}$, which is to say that the images of $a_{1}$ and $b_{1}$ in $\eta_{1,K}^{1}$ are equal. The diagram
  \begin{equation*}
    \begin{tikzcd}[column sep=huge]
      \eta_{1,K}^{0} \times \lambda_{1,K}^{0}
      \ar[r,
      "f_{2,K}^{\times, 0}"]
      \ar[d, swap,
      "d_{\eta_{1,K}^{\bullet}} \times d_{\lambda_{1,K}^{\bullet}}"]
      &
      \eta_{2,K}^{0} \times \lambda_{2,K}^{0}
      \ar[d,
      "d_{\eta_{2,K}^{\bullet}} \times d_{\lambda_{2,K}^{\bullet}}"]
      \\
      \eta_{1,K}^{1} \times \eta_{1,K}^{1}
      \ar[r, "f_{2,K}^{1}\times f_{2,K}^{1}"]
      &
      \eta_{2,K}^{1} \times \eta_{2,K}^{1}
    \end{tikzcd}
  \end{equation*}
  commutes because $f_{1,K}^{\bullet}$ and $\MP f_{1,K}^{\bullet}$ are morphisms of complexes. The image of $(a_{1}, b_{1})$ in the bottom right corner of this diagram has the form $(c,c)\in \eta_{2,K}^{1} \times \eta_{2,K}^{1}$ by the assumption that $(a_{2}, b_{2})$ defines a pair in the relevant fiber product. By the assumption that \eqref{eq:ses-of-complexes-2} is short exact, $f_{2,K}^{1}$ is injective. Since the diagram commutes, we see that images of $a_{1}$ and $b_{1}$  in $\eta_{1,K}^{1}$ are equal. This has shown that $(a_{1}, b_{1})$ satisfies $(a_{1}, b_{1})\in \eta^{0}_{1,K} \times_{\eta^{1}_{1,K}} \lambda^{0}_{1,K}$. Taking $x_{2}$ as above, we have
  \[
    f_{2, K}^{\times, 0}((a_{1}, b_{1}))
    = (a_{2},b_{2}) + \Delta^{-}_{2}(x_{2}),
  \]
  so that
  \[
    f_{2, K}^{+,0}([a_{1}, b_{1}])
    = [a_{2},b_{2}].
  \]
  Since $[a_{2}, b_{2}]$ was an arbitrary element of $\ker f_{3,K}^{+, 0}\subset (\eta_{2}^{\bullet}, \nu_{2})_{K}^{0}$,  we conclude that $\ker f_{3,K}^{+, 0}\subset \im f_{2, K}^{+, 0}$. By \Cref{coch-les-lemma}, which we noted before applies in the current situation (i.e., to the complex \eqref{eq:ses-of-complexes-2-K}), we conclude that the functor $(-)_{K}$ of \eqref{eq:gen-fiber-functor-def} is exact.
    \qed}
\end{proposition}

\begin{remark}
\label{N-equals-0-remark}
We see from \eqref{eq:nu-star-eta-0-explicit} that $\MP(\eta^{\bullet}, \nu)_{K}^{0} = \eta_{K}^{1}\oplus \eta_{K}^{-1}$  if $\nu = 0$. It follows that $(\eta^{\bullet},0)_{K}^{0} = \eta_{K}^{0}.$
\end{remark}

\begin{remark}
  Let $(\eta^{\bullet}, \nu)\in \ExtMon$. By \Cref{monodromy-pushout-is-1sst} and \Cref{baer-ht-closure}, $(\eta^{\bullet},\nu)_{K}^{0}$ is in $\Repsstonetors$.
\end{remark}

\subsection{Torsion in abelian $K$-varieties with semistable reduction}

The following result is proved by applying a more general result of Bertapelle, Candilera and Cristante to the specific context of \Cref{1-motive-gives-torsion}.

\begin{proposition}[\cite{bertapelle-candilera-cristante}, Theorem 19]
  \label{BCC-main-theorem}
Let $A_{K}$ be a semistable abelian $K$-variety, let $n\geq 2$ be a power of $p$. Let $\sQ_{K}$ be the log 1-motive of \Cref{1-motive-gives-torsion}, let $\sQ^{1}_{\varpi}$ be the  $\OK$-1-motive of \Cref{raynaud-decomp}, and let $\nu_{n}$ be the monodromy homomorphism \eqref{eq:nu-n-def}. There exists an isomorphism
  \[
    (\eta(\sQ^{1}_{\varpi}, n),\nu_{n})_{K} \simeq \eta(\sQ_{K},n)
  \]
  of objects of $\Ext_{\fin^{p}_{K}}(Y_{K}/nY_{K}, \widetilde A_{K}[n])$. In particular, the degree-zero term in the generic fiber $(\eta(\sQ^{1}_{\varpi}, n),\nu_{n})_{K}$ is isomorphic to $A_{K}[n]$.
\end{proposition}  

Let $(\eta^{\bullet}, \nu)\in \ExtMon$ and let $H = \eta^{0}$. Since $\eta^{1}$ is \'etale, there is a canonical morphism of complexes $\ceseq(H) \to \eta^{\bullet}$. In fact, as we shall see, $\nu$ induces a monodromy structure on $\ceseq(H)$. While it is not strictly necessary for us to work in $\fin^{p, N}_{\OK}$, it is desirable to do so because it allows us to frame our later results in terms of connected-\'etale sequences of objects, which are commonly used, and which are the context of Kato's work on log finite flat group schemes.

\begin{definition}
We call the functor $(-)_{K}$ of \Cref{generic-fiber-functor-is-exact} the generic fiber functor on $\ExtMon$. This defines an exact functor on $\fin^{p,N}_{\OK}$ by sending $(H,N)\in \fin^{p,N}_{\OK}$ to
  \[
(H, N)_{K} := (\ceseq(H), N)_{K}^{0}.
\]
\end{definition}

The following result lets us compare generic fibers of objects of $\ExtMon$ to generic fibers of objects of $\fin^{p, N}_{\OK}$.

\begin{proposition}
  \label{same-gen-fiber}
  Let $(\eta_{2}^{\bullet}, \nu_{2})\in \ExtMon$ and let  $\eta_{1}^{\bullet}$ be an acyclic cochain complex of objects of $\fin^{p}_{\OK}$ concentrated in degrees $-1$, $0$ and $1$, with $\eta_{1}^{-1}\subset \eta_{2}^{-1},$ $\eta_{1}^{0}=\eta_{2}^{0}$, and with $\eta_{2}^{1}$ \'etale. Then there exists a morphism $\nu_{1}$ making $(\eta_{1}^{\bullet}, \nu_{1})$ an object of $\ExtMon$, and the natural morphism $f^{\bullet}: \eta_{1}^{\bullet}\to \eta_{2}^{\bullet}$ is in fact a morphism $(\eta_{1}^{\bullet}, \nu_{1}) \to (\eta_{2}^{\bullet}, \nu_{2})$ inducing an isomorphism
  \[
    f^{+,0}_{K}: (\eta_{1}^{\bullet}, \nu_{1})_{K}^{0}
    \longisoto (\eta_{2}^{\bullet}, \nu_{2})_{K}^{0}.
  \]
\proof{
  We first define $\nu_{1}$. Recall that $\nu_{2}$ is a map $\nu_{2}:\eta_{2}^{1}(1)\to \eta_{2}^{-1}$. Since $\eta_{2}^{1}(1)$ is of multiplicative type and of $p$-power order, $\nu_{2}$ factors through the identity component $(\eta_{2}^{-1})^{\circ} = (\eta_{1}^{-1})^{\circ}$. Pre-composing $\nu_{2}$ with $f^{1}(1): \eta_{1}^{1}(1) \to \eta_{2}^{1}(1)$, we obtain a morphism
  \begin{equation}
    \label{eq:nu1-def}
    \nu_{1}: \eta_{1}^{1}(1)\to \eta_{1}^{-1}.
  \end{equation}
By construction of $\eta_{1}$, we obtain a morphism $f^{\bullet}: (\eta_{1}^{\bullet}, \nu_{1}) \to (\eta_{2}^{\bullet}, \nu_{2})$.

  We now show that $f_{K}^{+,0}$ is an isomorphism. We reuse our conventions and notation from the proof of \Cref{generic-fiber-functor-is-exact}, while also letting $E_{n} := E_{\varpi, K}[n]$.  We have a commutative diagram
  \begin{small}
\begin{equation}
  \label{eq:big-diagram}
  \begin{tikzcd}[column sep=tiny]
    0\ar[rr] 
    & & \eta_{1, K}^{1}(1) \ar[rr]   \ar[dr, two heads, "f_{K}^{1} \otimes \id"]
                                   \ar[dd, near start, "\nu_{1}"]
& & \eta_{1, K}^{1}\otimes E_{n}
                               \ar[rr] \ar[dr, two heads, "f_{K}^{1}\otimes \id"]
                               \ar[dd, near start, "\nu_{1}'"]
& & \eta_{1, K}^{1}              \ar[rr] \ar[dr, two heads, "f_{K}^{1}"]
                               \ar[dd, equals]
& & 0 
& \\
&   0 \ar[rr]
& &  \eta_{2, K}^{1}(1)   \ar[rr] \ar[dl, near start, "\nu_{2}"]
& & \eta_{2, K}^{1}\otimes E_{n}
                               \ar[rr] \ar[dd, near start, "\nu_{2}'"]
& & \eta_{2, K}^{1}\ar[rr]       \ar[dd, equals] 
& & 0  
\\
    0                         \ar[rr] 
& & \eta_{1, K}^{-1}             \ar[rr]\ar[dr, hook] 
& & \lambda_{1, K}^{0}                 \ar[rr]\ar[dr, "\MP f^{0}_{K}"] 
& & \eta_{1, K}^{1}             \ar[rr]\ar[dr, two heads] 
& & 0 
& \\
&   0\ar[rr] 
& &  \eta_{2, K}^{-1}           \ar[rr] 
& & \lambda_{2, K}^{0}               \ar[rr] 
& & \eta_{2, K}^{1}             \ar[rr] 
& & 0.
\end{tikzcd}
\end{equation}
\end{small}
The homomorphism $\MP f^{0}_{K}$ is induced by $\eta_{1, K}^{1}\otimes E_{n} \to \eta_{2, K}^{1}\otimes E_{n} \to \lambda_{2, K}^{0}$ and $\eta_{1, K}^{-1}\to \eta_{2, K}^{-1} \to \lambda_{2, K}^{0}$ via the universal property of $\lambda_{1, K}^{0}$ as a cofiber product. We note also that it need not be the case that $\MP f^{0}_{K}$ is an isomorphism, as can be seen in the example where $\eta_{2,K}^{1} = 0$.

We now show that $f_{K}^{+, 0}$ is an isomorphism of finite $K$-groups. For each $i\in \{1,2\}$, $\eta_{i, K}^{0}$ and $(\eta_{i}^{\bullet}, \nu_{i})_{K}^{0}$ both are objects of $\Ext_{\fin^{p}_{K}}(\eta_{i,K}^{1}, \eta_{i,K}^{-1})$, hence $\eta_{i}^{0}$ and $(\eta_{i}^{\bullet}, \nu_{i})_{K}^{0}$ each have the same order. Since $f^{0}$ is the identity map, we see that $(\eta_{1}^{\bullet}, \nu_{1})_{K}^{0}$ and $(\eta_{2}^{\bullet}, \nu_{2})_{K}^{0}$ have the same order. Therefore, to show that $f_{K}^{+,0}$ is an isomorphism, it suffices to show that it is an injection. Note that this does not follow from exactness of $(-)_{K}$ since it may be that $f^{1}$ is not injective.

Recall that any element of $(\eta_{1}^{\bullet}, \nu_{1})_{K}^{0}$ can be represented as $[a, [m, n]]$ with
\[
  (a, (m, n)) \in
  \eta_{1, K}^{0}
  \times_{\eta_{1, K}^{-1}} (\eta_{1, K}^{1}\otimes E_{n} \times \eta_{1, k}^{-1}).
\]
Then
$f_{K}^{+,0}([a, [m, n]])
= [f_{K}^{0}(a), g([m, n])]
= [a, [f_{K}^{1}\otimes \id(m), n]]$. Now let $[a, [m, n]] \in \ker f_{K}^{+, 0}$. Then there exist $c \in \eta_{2, K}^{-1}$ and $y \in \eta_{2, K}^{1}$ such that
\begin{align}
  \label{eq:ker-condition}
    (a, (f_{K}^{1}\otimes \id(m), n)) = (c, (y, -c - \nu_{2}(y)))
  \in
  \eta_{2, K}^{0}
  \times_{\eta_{2, K}^{-1}} (\eta_{2, K}^{1}\otimes E_{n} \times \eta_{2, k}^{-1}).
\end{align}
Since $n$ and $\nu_{2}(y)$ both are in $\eta_{1, K}^{-1}$ and $c = -n - \nu_{2}(y)$, we see that $c$ is in $\eta_{1, K}^{-1}$. Then $a = c$ maps to 0 in $\eta_{1, K}^{1}$. Now, $(a, [m, n])$ defines a point in a fiber product over $\eta_{1, K}^{-1}$, so we see that the surjection $\lambda_{1, K}^{0}\to \eta_{1, K}^{1}$ maps to $[m, n]$ to $0$, and from \eqref{eq:big-diagram} we see that $m\in \eta_{1, K}^{1}(1)$.

Since $m\in \eta_{1, K}^{1}(1)$, we have $f_{K}^{1}\otimes \id(m)\in \eta_{2, K}(1)$. By \eqref{eq:ker-condition} we have $f_{K}^{1}\otimes \id(m) = y$, so we see that from \eqref{eq:big-diagram} that $\nu_{2}(y) = \nu_{2}(f_{K}^{1}\otimes \id(m)) = \nu_{1}(m).$ Using the relationships derived above, we have shown that there exists $c\in \eta_{1, K}^{-1}$ such that, with $\Delta_{1}^{-}$ defined as in the proof of \Cref{generic-fiber-functor-is-exact}, the equation 
\begin{align*}
  [a, [m, n]]
  = [c, [m, -c-\nu_{2}(y)]] 
  = \Delta_{1}^{-}(c).
\end{align*}
holds. This has shown that $[a, [m, n]] = 0$.
\qed}\end{proposition}

\begin{corollary}
  There exists a canonical, essentially surjective functor $\ExtMon\to \fin^{p, N}_{\OK}$ that sends an object $(\eta^{\bullet}, \nu)$ to $(\eta^{0}, \nu|_{\eta^{0,\et}(1)})$. The natural map $(\eta^{\bullet}, \nu)_{K}^{0} \to (\eta^{0}, \nu|_{\eta^{0,\et}(1)})_{K}^{0} $ is an isomorphism.
  \proof{This is the special case of \Cref{same-gen-fiber} where $\eta_{1}^{\bullet} = \ceseq(\eta_{2}^{0})$.    
    \qed}
\end{corollary}

\begin{corollary}
  \label{BCC-N-cor}
Let $A_{K}$ be a semistable abelian $K$-variety, let $n\geq 2$ be a power of $p$. Let $\sQ_{K}$ be the log 1-motive of \Cref{1-motive-gives-torsion}, let $\sQ^{1}_{\varpi}$ be the  $\OK$-1-motive of \Cref{raynaud-decomp}. Attached to our uniformizer $\varpi$ of $\OK$ there is an object $(\sQ^{1}_{\varpi}[n], N_{n})$ of $\fin^{p, N}_{\OK}$ with a canonical isomorphism
  \[(\sQ^{1}_{\varpi}[n], N_{n})_{K} \simeq A_{K}[n].\]
\end{corollary}

\subsection{Monodromy criterion for torsion 1-crystallinity}

Let $(Q,N)$ be an object of $\finpN$, and let $V = (Q,N)_{K} \in \Reponessttors$. Informed by Fontaine's theory of filtered $(\varphi,N)$-modules, it is natural to seek a monodromy criterion for 1-crystallinity of subobjects of $V$. One might hope to prove such a criterion by first developing an analogous theory of monodromy structures on integral objects, and then showing that $(Q,N)$, comes from such an object. While \cite{kisin-f-crystals} does develop a theory of Kisin modules with monodromy, the monodromy operators there are only defined after inverting $p$. In an unpublished preprint \cite{madapusi-pera-log-pdivs-and-sst-reps}, Madapusi Pera has developed a theory of height-1 Kisin modules with integral monodromy and has shown that every object of $\Reponesst$ is the Tate module of some log $p$-divisible group. It seems possible that one could then prove a strengthening of Raynaud's embedding theorem \cite{raynaud-pp} 2.3.1 for log finite flat $\OK$-group schemes, and then to apply results from Fontaine's theory to prove a monodromy criterion. We take a different approach by working directly with objects of $\ExtMon$ and applying the Torsion Full Faithfulness Theorem \Cref{torsion-ff}.

Fix $(\eta^{\bullet}, \nu) \in \ExtMon$ until \Cref{sec:main-results}.

  \begin{proposition}
  \label{1-crys-iff-1-i-star-n-star-crys}
The $\Gamma_{K}$-module $(\eta^{\bullet}, \nu)_{K}^{0}$ is 1-crystalline if and only if $\MP(\eta^{\bullet}, \nu)_{K}^{0}$ is 1-crystalline.
  \proof{
This is an application of \Cref{baer-ht-closure}. \qed}
\end{proposition}

Let $U^{\bullet} = (\eta^{\bullet}, \nu)_{K}$ and let $\nu_{*}U = \MP((\eta^{\bullet}, \nu))_{K}^{0}$. Recall that $\nu_{*}U$ fits into a commuting diagram
\begin{equation}
  \label{eq:push-pull}
    \begin{tikzcd}
    0 \ar[r]
    & U^{1}(1)\ar[r] \ar[d, "\nu"]
    & E_{\varpi, K}[n]\otimes U^{1} \ar[r] \ar[d, "\nu"]
    & U^{1}\ar[r] \ar[d, equals]
    & 0
    \\
    0 \ar[r]
    & U^{-1} \ar[r]
    & \nu_{*}U \ar[r]
    & U^{1} \ar[r]
    & 0
  \end{tikzcd}
\end{equation}
where the leftmost square is cocartesian, where $n$ is the order of $\eta^{1}$, and where the top row is obtained by tensoring $U^{1}$ against the short exact sequence
\[
  \tikzses{\theta^{\varpi}_{n,K}:}{\mu_{n,K}}{}{E_{\varpi, K}[n]}{}{\Z/n\Z_{K}}.
\]
Choose an ordered basis $\{x_{1}, x_{2}\}$ of $E_{\varpi, K}[n]$ as in \Cref{subsec:q-tate-eg} (where now we have specified that $q=\varpi$). In particular, $x_{1}$ spans the $\Gamma_{K}$-stable line in $E_{\varpi, K}[n]$. If $n=p^{m}$, the action of $\Gamma_{K}$ on $x_{2}$ is determined by a $\Z/n\Z$-valued cocycle $c_{\varpi, n}$ defined by the relation
\begin{align}
  \label{eq:c-varpi-n-relation}
  \sigma(\varpi^{\flat, (m)}) = (\epsilon^{(m)})^{c_{\varpi, n}(\sigma)}\varpi^{\flat, (m)},
\end{align}
for all $\sigma\in \Gamma_{K}$, where $\varpi^{\flat}$ and $\epsilon$ are as in \Cref{subsec:notation}. We record the following basic lemma.
\begin{lemma}
  \label{p-nmid-cocycle}
  There exists $\sigma\in \Gamma_{K}$ such that $c_{\varpi, n}(\sigma)\in (\Z/n\Z)^{\times}$.
  \proof{
    Suppose $p \mid c_{\varpi, n}(\sigma)$ for all $\sigma$. Recall that $n=p^{m}$ with $m\geq 1$. Raising the relation defining $c_{\varpi, n}$ to the $p^{m-1}$-th power, we see that $(\varpi^{\flat, (m)})^{p^{m-1}}$ is fixed by all of $\Gamma_{K}$, which is to say that $K$ contains a $p$-th root of $\varpi$. Since $\varpi$ is a uniformizer of $\OK$, this is impossible, so we conclude that $c_{\varpi, n}(\sigma)\in (\Z/n\Z)^{\times}$ for some
$\sigma$.
  \qed}
\end{lemma}

With our choice of basis of $E_{\varpi, K}[n]$, we may write elements of the tensor product
$E_{\varpi, K}[n]\otimes U^{1}$ uniquely in the form
\[
  x_{1}\otimes y_{1} + x_{2}\otimes y_{2} \qquad (y_{1},y_{2}\in U^{1}).
\]
Using this basis for $E_{\varpi, K}[n]$, the two morphisms in the top row of \eqref{eq:push-pull} are simply the inclusion
\begin{align*}
  \twoar{}{
  U^{1}(1) = \Z/n\Z(1) \otimes U^{1}}{
  E_{\varpi, K}[n]\otimes U^{1}}{x_{1}\otimes y_{1}}{x_{1}\otimes y_{1}}
\end{align*}
and the morphism
\begin{align*}
  \twoar{}{E_{\varpi, K}[n]\otimes U^{1}}{U^{1}}{x_{1}\otimes y_{1}+x_{2}\otimes y_{2}}{y_{2}.}
\end{align*}

As in \eqref{eq:eta-i-baer-desc}, $\nu_{*}U$ satisfies
\begin{align}
  \label{eq:N-star-W-baer-def}
\nu_{*}U = \frac{(E_{\varpi, K}[n]\otimes U^{1})\oplus U^{-1}}{\{(x_{1}\otimes y_{1}, 0)-(0, \nu(x_{1}\otimes y_{1})):\: y_{1}\in U^{1} \}}.  
\end{align}
Given $(x_{1}\otimes y_{1} + x_{2}\otimes y_{2}, z)\in (E_{\varpi, K}[n]\otimes U^{1})\oplus U^{-1}$, we denote its class in $\nu_{*}U$ by $[x_{1}\otimes y_{1} + x_{2}\otimes y_{2}, z]$. Using the relation in \eqref{eq:N-star-W-baer-def}), we have
\[
\nu_{*}U = \{[x_{2}\otimes y_{2}, z]:\: y_{2}\in U^{1},\, z\in U^{-1}\}.
\]
There is a homomorphism of finite abelian groups
\begin{align}
  \label{eq:s-def} \twoardef{s}{U^{1}}{\nu_{*}U}{u}{
  [x_{2}\otimes u,0]}
\end{align}
which is a section of the quotient map $\nu_{*}U\to U^{1}$. The morphism $s$ need not be $\Gamma_{K}$-equivariant, but we have the following result.
\begin{lemma}
  The homomorphism $s$ is $\Gamma_{\varpi^{\flat}}$-equivariant.
  \proof{Restrict \eqref{eq:push-pull} to $\Gamma_{\varpi^{\flat}}$. From the relation \eqref{eq:c-varpi-n-relation} defining $c_{\varpi, n}$, we see that $c_{\varpi, n}|_{\Gamma_{\varpi^{\flat}}} = 0$, so that our choice of basis $\{x_{1}, x_{2}\}$ induces an isomorphism
\begin{align*}
  \Resoo(E_{\varpi, K}[n])
  \simeq \Resoo\big(\Z/n\Z(1)\oplus \Z/n\Z\big)
\end{align*}
We see that the homomorphism $y_{2}\mapsto x_{2}\otimes y_{2}$ gives a $\Gamma_{\varpi^{\flat}}$-equivariant splitting of the top row of \eqref{eq:push-pull}. Composing with $\nu$, we obtain the map $s$, which we see is a $\Gamma_{\varpi^{\flat}}$-equivariant splitting of
\begin{align*}
 \ses{\Resoo U^{-1}}{}{\Resoo \nu_{*}U}{}{\Resoo U^{1}}.
\end{align*}
    \qed}
\end{lemma}

\begin{proposition}
  \label{1-crys-iff-equivariant}
  The $\Gamma_{K}$-module $\nu_{*}U$ is 1-crystalline if and only if the $\Gamma_{\varpi^{\flat}}$-equivariant section $s$ defined in \eqref{eq:s-def} is $\Gamma_{K}$-equivariant.
  \proof{Recall that a torsion $\Gamma_{K}$-module is 1-crystalline if and only if it is the generic fiber of a finite flat $\OK$-group. Therefore $U^{-1} = \eta_{1,K}$ and $U^{1} = \eta_{-1,K}$ are 1-crystalline. Suppose first that $s$ is $\Gamma_{K}$-equivariant. Then $\nu_{*}U \simeq U^{-1}\oplus U^{1}$ as $\Zp[\Gamma_{K}]$-modules. Since finite direct sums of 1-crystalline $\Gamma_{K}$-modules are 1-crystalline, we conclude that $\nu_{*}U$ is 1-crystalline.

    Now suppose $\nu_{*}U$ is 1-crystalline. Since $U^{1}$ is 1-crystalline, \Cref{torsion-ff} implies that the $\Gamma_{\varpi^{\flat}}$-equivariant morphism $s$ is $\Gamma_{K}$-equivariant. 
  \qed}
\end{proposition}

\begin{remark}
  In \Cref{1-crys-iff-equivariant} we made non-canonical choices. On the one hand, $\varpi^{\flat}$ gives the subgroup $\Gamma_{\varpi^{\flat}}\subset \Gamma_{K}$ that was used when applying \Cref{torsion-ff}. On the other hand, $\varpi = \varpi^{\flat, (0)}$ was used when we formed $E_{\varpi, K}[n]$ and performed Raynaud decomposition. We could have chosen a different uniformizer, say $u\varpi$ where $u\in \OK^{\times}$, to perform Raynaud decomposition. In that case, $E_{u\varpi, K}[n]$ need not be split as a $\Zp[\Gamma_{\varpi^{\flat}}]$-module.
\end{remark}

Returning to the situation of the above \Cref{1-crys-iff-equivariant}, we now relate the splitness of $\nu_{*}U$ as an extension of $U^{1}$ by $U^{-1}$ to the monodromy morphism. Recall that $\nu$ gives a morphism $\nu(-1):U^{1}\to U^{-1}(-1)$ of $\Gamma_{K}$-modules.

\begin{proposition}
  \label{equivariant-iff-in-ker}
  The homomorphism $s$ is $\Gamma_{K}$-equivariant if and only if $\nu(-1)=0$.
  \proof{
    Suppose $\nu(-1) = 0$. In the basis $\{x_{1}, x_{2}\}$ for $E_{\varpi, K}[n]\otimes U^{1}$, this means that $\nu(x_{1}\otimes u)=0$ for all $u\in U^{1}$. Let $\sigma\in \Gamma_{K}$. We have 
    \begin{align*}
      \sigma s(u) - s(\sigma u)
      &= \sigma [x_{2}\otimes u,0]
        - [x_{2}\otimes \sigma u,0] \\
      &= [c_{\varpi, n}(\sigma)x_{1}\otimes \sigma u, 0] \\
      &= [0, \nu(x_{1}\otimes c_{\varpi, n}(\sigma)\sigma u)] 
      = [0, 0].
    \end{align*}
    This shows that $s$ is $\Gamma_{K}$-equivariant.

    Conversely, suppose $s$ is $\Gamma_{K}$-equivariant and let $u\in U^{1}$. For any $\sigma\in \Gamma_{K}$, we see from the above calculation that we have an equality of elements of $\nu_{*}U$
    \[[0,0] = \sigma s(u) - s(\sigma u) = [0, \nu(x_{1}\otimes c_{\varpi, n}(\sigma)\sigma u)],\]
    which is to say that, for some $y_{1}\in U^{1}$ depending on $u$, we have the equality of elements of $(E_{\varpi, K}[n]\otimes U^{1})\oplus U^{-1}$
    \[
      (0, \nu(x_{1}\otimes c_{\varpi, n}(\sigma)\sigma u))
      = (x_{1}\otimes y_{1}, -\nu(x_{1}\otimes y_{1})).
    \]
    We see that we must have $y_{1}=0$, hence $\nu(x_{1}\otimes y_{1})=0$, hence
    \begin{align}
      \label{eq:N-kills-something}
\nu(x_{1}\otimes c_{\varpi, n}(\sigma)\sigma u) = 0.      
    \end{align}
    By \Cref{p-nmid-cocycle}, we may choose $\sigma$ such that $c_{\varpi, n}(\sigma)\in (\Z/n\Z)^{\times}$. Let $\chi$ be the $p$-adic cyclotomic character on $\Gamma_{K}$. It follows from \eqref{eq:N-kills-something} that
    \[
      0 = \nu(x_{1}\otimes \sigma u)
      = \chi(\sigma^{-1})\sigma \nu(x_{1}\otimes u),
    \]
whence $\nu(x_{1}\otimes u)=0$. We conclude that $\nu(-1) = 0$.
  \qed}
\end{proposition}

\begin{proposition}
  \label{monodromy-criterion}
  The $\Gamma_{K}$-module $(\eta^{\bullet}, \nu)_{K}^{0}$ is 1-crystalline if and only if $\nu = 0$.
  \proof{This is proved by combining \Cref{equivariant-iff-in-ker}, \Cref{1-crys-iff-equivariant}, and \Cref{1-crys-iff-1-i-star-n-star-crys}.
    \qed}
\end{proposition}

Let $(Q,N)\in \fin^{p, N}_{\OK}$. We now prolong $\Crys_{1}((Q, N)_{K})$ in $Q$.

\begin{proposition}
  \label{pullback-prolong-lemma}
  Let $(Q, N)$ be an object of $\finpN$ and let $V^{\bullet} = (\ceseq(Q), N)_{K}$. Let $U^{0}\subset V^{0} = (Q,N)_{K}$ be a $\Gamma_{K}$-submodule containing $V^{-1} = Q^{\circ}_{K}$. Let $\eta^{1}$ be the unique \'etale prolongation of $U^{0}/V^{-1}$, let $i^{1}:\eta^{1}\to Q^{\et}$ be the unique prolongation of the inclusion $U^{0}/V^{-1}\subset V^{1}$, and let $\eta^{\bullet} := (i^{1})^{*}\ceseq(Q)$. Then $\eta^{0}$ is the scheme theoretic closure of $U^{0}$ in $Q$, and there is a natural morphism $i^{\bullet}:\eta^{\bullet}\to \ceseq(Q)$ with the property that $i^{0}$ is an open immersion. Moreover, $(\eta^{\bullet}, N|_{\eta^{1}(1)})$ is a subobject of $(\ceseq(Q), N)$ via $i^{\bullet}$, and the morphism $i_{K}^{0}$ identifies $(\eta^{\bullet}, N|_{\eta^{1}(1)})_{K}^{0}$ with $U^{0}\subset V^{0}$.
  \proof{
    Since $\ceseq(Q)^{1} = Q^{\et}$ is \'etale, $V^{1} = Q^{\et}_{K}$ is unramified, so existence and uniqueness of $i^{1}$ are standard. By definition, $\eta^{0} = Q \times_{Q^{\et}} \eta^{1}$, and $i^{0}$ is simply the projection onto the first factor. Since $i^{1}$ is an open immersion, its base change $i^{0}$ is an open immersion. It is clear that $\eta^{0}$, the union of the components of $Q$ containing the points of the finite \'etale $K$-group $U^{0}$, is the scheme theoretic closure of $U^{0}$ in $Q$. Being that $\eta^{0}$ a union of components of $Q$, restricting $N$ to $\eta^{1}(1)$ (cf.~the construction of \eqref{eq:nu1-def}) gives a subobject $(\eta^{\bullet}, N|_{\eta^{1}(1)})$ of $(\ceseq(Q), N)$. By exactness of the generic fiber functor on $\finpN$ (\Cref{generic-fiber-functor-is-exact}), $i^{0}$ injects $(\eta^{\bullet}, N|_{\eta^{1}(1)})_{K}^{0}$ onto a $\Gamma_{K}$-submodule of $(\ceseq(Q), N)_{K}^{0} = (Q, N)_{K} = V^{0}$. By construction, $\im i_{K}^{0}$ and $U^{0}$ each consist of the same cosets of $V^{-1} = Q^{\circ}_{K}$ in $V^{0}$. 
    \qed}
\end{proposition}

\begin{proposition}
  \label{star-Q-realizes-crys-1}
  Let $(Q, N)$ be an object of $\finpN$. Let ${}^{*}\ceseq(Q)$, be the pullback of $\ceseq(Q)$ under the inclusion $\ker N(-1)\subset Q^{\et}$, so that ${}^{*}\ceseq(Q)$ is a class in $\Ext_{\finp}(\ker N(-1), Q^{-1})$. Then the natural map ${}^{*}\ceseq(Q)\to \ceseq(Q)$ gives rise to an isomorphism of $\Gamma_{K}$-modules $({}^{*}\ceseq(Q))_{K}^{0}\simeq \Crys_{1}((Q,N)_{K})$.
  \proof{By \Cref{N-equals-0-remark}, $({}^{*}\ceseq(Q), N|_{({}^{*}\ceseq(Q))^{1}})_{K} = ({}^{*}\ceseq(Q))_{K}$. As the $({}^{*}\ceseq(Q))^{0}$ is finite flat, exactness of the generic fiber functor implies that the submodule $({}^{*}\ceseq(Q))_{K}$ of $(Q, N)_{K}$  is contained in $\Crys_{1}((Q,N)_{K})$. We show the reverse inclusion. Since $Q^{\circ}$ is a finite flat $\OK$-group, we have $Q^{\circ}_{K} \subset \Crys_{1}((Q,N)_{K})$. Apply \Cref{pullback-prolong-lemma} to obtain a subobject $(\eta^{\bullet}, N|_{\eta^{1}(1)})$ of $(\ceseq(Q), N)$. By \Cref{monodromy-criterion}, $\eta^{1}(1)\subset \ker N$, so $\eta^{\bullet}$ is a subobject of ${}^{*}\ceseq(Q)$. Taking generic fibers, we have $\Crys_{1}((Q,N)_{K}) = \eta^{0}_{K}\subset {}^{*}\ceseq(Q)_{K}$. This completes the proof.
    \qed}
\end{proposition}

\section{Results on semistable abelian varieties}
\label{sec:main-results}

For this section, we adopt the notation of \Cref{subsec:notation}. In particular,  $p$ is any
rational prime and $K/\Qp$ a discretely valued extension. Let $A_{K}$ be an abelian $K$-variety with semistable reduction, let $A$ be the N\'eron model of $A_{K}$ and let $\widetilde A$ be the Raynaud extension.

\subsection{Realizing torsion 1-crystallinity, and the Coleman-Iovita criterion.}
Let $\sQ_{K}$ be the log 1-motive attached to $A_{K}$ using the functor $F(\OK,\fm_{\OK})$, as in \Cref{1-motive-gives-torsion}. Let  $n\geq 2$ be a power of $p$, and let $(\sQon, N_{n})\in \ExtMon$ and $(\sQon, N_{n}) \in \finpN$ be as in \Cref{BCC-main-theorem} and \Cref{BCC-N-cor}, respectively.  As we did in \Cref{star-Q-realizes-crys-1}, we define
\[
{}^{*}\sQon := 
  \sQon
  \times_{\sQon^{\et}}
  (\ker N_{n})(-1).
\]
Since $Y/nY$ is \'etale, passing to connected parts in $\eta(\sQ^{1}_{\varpi}[n],n)$ (an object of the group $ \Ext_{\fin^{p}_{\OK}}(Y/nY, \widetilde A[n])$) shows that the morphism $\widetilde A[n]\to \sQo[n]$ is an isomorphism on identity components, so that ${}^{*}\sQon$ fits into a diagram
\begin{equation}
  \label{eq:star-Q-1-varpi-pullback-diagram}
    \begin{tikzcd}
    0\ar[r] & \widetilde A[n]^{\circ}  \ar[r]\ar[d, equal] & {}^{*}\sQon \ar[r] \ar[d, hook] & (\ker N_{n})(-1) \ar[r]\ar[d, hook] & 0 \\
    0\ar[r] & \widetilde A[n]^{\circ} \ar[r] & \sQon \ar[r] & \sQon^{\et} \ar[r] & 0
  \end{tikzcd}
\end{equation}
where the rows are connected-\'etale sequences and the vertical morphisms are inclusions of open $\OK$-subgroups.
\begin{theorem}
  \label{main-theorem-torsion}
  The canonical morphism $({}^{*}\sQon, 0)\to (\sQon, N_{n})$ coming from \eqref{eq:star-Q-1-varpi-pullback-diagram} induces an isomorphism of $\Gamma_{K}$-modules ${}^{*}\sQon_{K}\simeq \Crys_{1}(A_{K}[n])$.
  \proof{
This follows directly from \Cref{star-Q-realizes-crys-1}.
    \qed}\end{theorem}
 \Cref{main-theorem-torsion} gives rise to a description of $\Crys_{1}(T_{p}A_{K})$. We note that this result about $T_{p}A_{K}$ follows from \cite{coleman-iovita-1}, Corollary 4.6, which is a key ingredient in their proof of their monodromy criterion for good reduction.

\begin{theorem}
  \label{main-theorem-free}
The isomorphisms  ${}^{*}\sQ^{1}_{\varpi}[p^{m}]_{K}  \simeq \Crys_{1}(A_{K}[p^{m}])$ for varying $m\in \Z_{\geq 1}$ are compatible in $m$, and give rise to an isomorphism
  \[T_{p}\widetilde A_{K} \simeq \Crys_{1}(T_{p}A_{K}).\]
  \proof{Recall that, for $m\in \Z_{\geq 1}$, the morphism $N_{p^{m}}: \sQ^{1}_{\varpi}[p^{m}]^{\et}(1) \to \sQ^{1}_{\varpi}[p^{m}]^{\circ}$ is obtained by pre-composing the morphism $\nu_{p^{m}}$ of \eqref{eq:nu-n-def} with the canonical morphism $\sQ^{1}_{\varpi}[p^{m}]^{\et}(1)\to (Y/p^{m}Y)(1)$. It is clear that $(\ker N_{p^{m}})_{m\geq 1}$ is an inverse system of objects of $\fin^{p}_{\OK}$. The isomorphisms of \Cref{main-theorem-torsion} are compatible for varying $p^{m}$; each such isomorphism is the composition of the projection $\pr_{1,K}$ on $\sQ_{K}[p^{m}] \times_{\sQ^{1}_{\varpi}[p^{m}]^{\et}_{K}} \ker N_{p^{m},K}(-1)$ with the isomorphism $\sQ_{K}[p^{m}]\simeq A_{K}[p^{m}]$ of \Cref{1-motive-gives-torsion}, and the isomorphisms of \Cref{1-motive-gives-torsion} can be seen from the proof to be compatible. By \Cref{crys-r-functor},
\[\Crys_{1}(T_{p}A_{K})
  = \lim_{m}\Crys_{1}(T_{p}A_{K}/p^{m}T_{p}A_{K})
  = \lim_{m}\Crys_{1}(A_{K}[p^{m}]),
\]
so by \Cref{main-theorem-torsion} we have isomorphisms
\begin{align}
  \label{eq:crys-1-TpAK-is-lim}
  \Crys_{1}(T_{p}A_{K})
  \simeq \lim_{m} {}^{*}\sQo[p^{m}]_{K}.  
\end{align}
    Since $T_{p}\widetilde A[p^{\oo}]$ is 1-crystalline and is a submodule of $T_{p}A_{K}$ (via the sequences $\eta(A_{K},n)$ of \Cref{1-motive-gives-torsion}), it is contained in $\Crys_{1}(T_{p}A_{K})$. The quotient $\Crys_{1}(T_{p}A_{K})/p^{m}\Crys_{1}(T_{p}A_{K})$ is torsion 1-crystalline, as is $\widetilde A[p^{m}]_{K}$, so we find that
\[
      \widetilde A[p^{m}]_{K}
      \subset \Crys_{1}(T_{p}A_{K})/p^{m}\Crys_{1}(T_{p}A_{K})
      \subset \Crys_{1}(A[p^{m}]_{K}) = {}^{*}\sQo[p^{m}]_{K}.
\]
The upper row in \Cref{eq:star-Q-1-varpi-pullback-diagram} is evidently the connected-\'etale sequence of ${}^{*}\sQ[n]$. It follows from \eqref{eq:R-1-motive-et-ses} that ${}^{*}\sQ[p^{m}]^{\et} = \ker N_{p^{m}}(-1)$ fits into the short exact sequence
\begin{align*}
  \ses{\widetilde A[p^{m}]^{\et}}{}{{}^{*}\sQo[p^{m}]^{\et}}{\pi_{n}}{\ker \nu_{p^{m}}(-1)}.
\end{align*}
By \Cref{ker-gives-component-group} (since $\Phi_{A_{K}}$ is finite), $T_{p}((\ker \nu_{p^{m}}(-1))_{m\geq 1})=0$. It follows that $\lim_{m}{}^{*}\sQo[p^{m}]_{K} = \lim_{m}\widetilde A_{K}[p^{m}] = T_{p}\widetilde A_{K}$. By \eqref{eq:crys-1-TpAK-is-lim}, we have $T_{p}\widetilde A_{K}\simeq \Crys_{1}(T_{p}A_{K})$. This completes the proof.
\qed}
\end{theorem}

\begin{corollary}[\cite{coleman-iovita-1} 4.7; also \cite{breuil-modules-filtrees} 5.3.4 when $p>2$ (see also Breuil's Remarque 5.3.5 for earlier special cases)]
  \label{coleman-iovita-breuil}
  An abelian $K$-variety $A_{K}$ has good reduction over $\OK$ if and only if its $p$-adic Tate module $T_{p}A_{K}$ is 1-crystalline.
  \proof{
Our proof essentially follows the proof of Theorem 4.7 of \cite{coleman-iovita-1}. The primary difference is that we will invoke \Cref{main-theorem-free} where they invoke their Corollary 4.6, which they prove in quite different terms.

    Suppose that $A_{K}$ has good reduction. Then $T_{p}A_{K}$ is the generic fiber of the $p$-divisible group $A[p^{\oo}]$ over $\OK$, which is 1-crystalline by \Cref{1-crys-things-are-gen-fibers}.

For the converse, suppose that $T_{p}A_{K}$ is 1-crystalline. In the proof of Theorem 4.7 in \cite{coleman-iovita-1}, it is shown using the results 7.3.3 and 7.5.1 of \cite{fontaine-modules-galoisiens} that we may assume without loss of generality that $A_{K}$ has semistable reduction. Make this assumption.  By \Cref{torsion-in-1-motive} and \Cref{1-motive-gives-torsion}, there is a short exact sequence
  \[
    \ses{T_{p}\widetilde A_{K}}{}{T_{p}A_{K}}{}{Y_{K}\otimes \Qp/\Zp},
  \]
  where $Y_{K}$ has $\Zp$-rank equal to the dimension of the torus in $\widetilde A_{K}$. By \Cref{main-theorem-free}, the first morphism in this short exact sequence is an isomorphism. Therefore $Y = 0$, so $T = 0$ and $\widetilde A = B$. Since $\widetilde A_{K} = A_{K}$, we see that $B$, which is an abelian $\OK$-scheme, has generic fiber $A_{K}$. By \cite{blr} 1.2/8, $B$ is the N\'eron model of $A_{K}$. Thus we see that $A_{K}$ has good reduction if $T_{p}A_{K}$ is 1-crystalline.
    \qed}
\end{corollary}

\subsection{Results on N\'eron component groups}
Recall that $\Phi_{A_{K}}$ is the N\'eron component group of $A_{K}$, an \'etale $\OK$-group (\Cref{def:connected-neron}). It is shown in \cite{grothendieck-monodromy} that, using $(-)^{\fpart}$ for the finite part of a quasi-finite flat $\OK$-group, one has isomorphisms
\begin{align}
  \label{eq:phi-p-m-ito-fparts}
    \Phi_{A_{K}}[l^{m}] \simeq \frac{(A[l^{m}]^{\fpart})_{K}}{(A^{\circ}[l^{m}]^{\fpart})_{K}} \qquad (l\neq p,\, m\geq 1).
\end{align}
Kim and Marshall (\cite{kim-marshall}, p. 611) noted that this formula holds also when $l=p$, provided one assumes semistable reduction. In the case where $p>2$ and $K$ is unramified, they use this to prove \eqref{eq:Phi-pn-ito-crys} below. We give a different proof\footnote{To compare with \eqref{eq:phi-p-m-ito-fparts}, one uses the fact (\cite{lan-PUP}, 3.4.2.1) that $\widetilde A[p^{m}] \simeq A[p^{m}]^{\fpart} $.} of \eqref{eq:Phi-pn-ito-crys} using \Cref{main-theorem-torsion} and \Cref{main-theorem-free}. Let $n\geq 2$ be a power of $p$.

\begin{lemma}
  \label{main-lemma-crys-quos}
  There is a canonical isomorphism of unramified finite $\Gamma_{K}$-modules
\begin{align}
  \label{eq:Phi-pn-ito-crys}
  \Phi_{A_{K}}[n] \simeq
  \frac{\Crys_{1}(T_{p}A_{K}\otimes \Z/n\Z)}{ \Crys_{1}(T_{p}A_{K})\otimes \Z/n\Z}.
\end{align}
\proof{
Recall that ${}^{*}\sQon$ is defined by pullback as in
\[
  \begin{tikzcd}
    0\ar[r] & \widetilde A[n]^{\circ}  \ar[r]\ar[d, equal] & {}^{*}\sQon \ar[r] \ar[d] & (\ker N_{n})(-1) \ar[r]\ar[d,hook] & 0 \\
    0\ar[r] & \widetilde A[n]^{\circ} \ar[r] & \sQon \ar[r] & \sQon^{\et} \ar[r] & 0.
  \end{tikzcd}
\]
By the universal property of ${}^{*}\sQon$ as a fiber product, we obtain a homomorphism $\widetilde
A[n]\to {}^{*}\sQon$. We claim that the generic fiber of the morphism $\widetilde A[n]\to {}^{*}\sQon$ has cokernel isomorphic to $\Phi_{A_{K}}[n]$. Recall the short exact sequence
\[
  \ses{\widetilde A[n]^{\et}}{}{\sQon^{\et}}{\pi_{n}}{Y/nY}
\]
of \eqref{eq:R-1-motive-et-ses}, and recall that $N_{n}$, defined as in \Cref{same-gen-fiber}, factors as
\[
    N_{n}: \sQon^{\et}(1)\overset{\pi_{n}(1)}\longto Y/nY(1)\overset{\nu_{n}}\longto T[n] \hooklongrightarrow \sQon^{\circ}.
  \]
  By restricting $\pi_{n}(1)$ and untwisting, we obtain a short exact sequence
  \[
    \  \ses{\widetilde{A} [n]^{\et}}{}{(\ker N_{n})(-1)}{
    }{\ker \nu_{n}(-1)},
  \]
where the morphism $\widetilde{A} [n]^{\et} \to (\ker N_{n})(-1)$ is the one obtained by descending the canonical homomorphism $\widetilde A[n]\to {}^{*}\sQon$ to \'etale parts. We see that

\[
{}^{*}\sQon \longto (\ker N_{n})(-1) \longto \frac{(\ker N_{n})(-1)}{\widetilde A[n]^{\et}}
\]
induces an isomorphism
\begin{align}
  \label{eq:star-sQon-wtA-cpt-gp-intmdt}
    \frac{{}^{*}\sQon}{\widetilde A[n]}
  \simeq \frac{(\ker N_{n})(-1)}{\widetilde A[n]^{\et}}
  = \ker \nu_{n}(-1).
\end{align}
Via the isomorphism $(\sQ^{1}_{\varpi}[n], N_{n})^{0}_{K}\simeq A_{K}[n]$ of \Cref{BCC-N-cor}, we obtain isomorphisms
\[
  {}^{*}\sQon \simeq \Crys_{1}(A_{K}[n])
  \quad \text{ and } \quad
  T_{p}\tilde A_{K} \simeq \Crys_{1}(T_{p}A_{K})
\]
as in \Cref{main-theorem-torsion} and \Cref{main-theorem-free}, respectively. Using these, we obtain from \eqref{eq:star-sQon-wtA-cpt-gp-intmdt} and the canonical isomorphism $\ker \nu_{n}(-1) \simeq \Phi_{A_{K}}[n]$ of \Cref{ker-gives-component-group} an isomorphism
\[  \Phi_{A_{K}}[n] \simeq
  \frac{\Crys_{1}(T_{p}A_{K}\otimes \bZ/n\bZ)}{ \Crys_{1}(T_{p}A_{K})\otimes \bZ/n\bZ}.\]
This completes the proof
\qed} \end{lemma}

\begin{theorem}
  \label{main-theorem-component-groups}
  Let $\Phi_{A_{K}}[p^{\oo}]$ denote the $p$-Sylow subgroup of $\Phi_{A_{K}}$. Then the isomorphisms of \Cref{main-lemma-crys-quos} give rise to an isomorphism
  \begin{align}
  \label{eq:R1Crys1-torsion}
  \Phi_{A_{K}}[p^{\oo}] \simeq (R^{1}\Crys_{1}(T_{p}A_{K}))_{\tors}.
\end{align}
\proof{
  By \Cref{crys-and-colim},
  \[
    \colim_{m\geq 1} \Crys_{1}(T_{p}A_{K}\otimes \Z/p^{m}\Z)
    = \Crys_{1}(T_{p}A_{K}\otimes \Qp/\Zp).
  \]
  Hence, taking colimits, \eqref{eq:Phi-pn-ito-crys} gives
\begin{align}
  \label{eq:Phi-ito-crys}
  \Phi_{A_{K}}[p^{\oo}] =
  \frac{\Crys_{1}(T_{p}A_{K}\otimes \Qp/\Zp)}{ \Crys_{1}(T_{p}A_{K})\otimes \Qp/\Zp}.
\end{align}
On the other hand, consider the long exact sequence
\begin{equation*}
  \begin{tikzcd}[column sep=tiny]
    0 \ar[r]
  & \Crys_{1}(T_{p}A_{K}) \ar[r]
  & \Crys_{1}(T_{p}A_{K})\otimes \Qp \ar[r]
  & \Crys_{1}(T_{p}A_{K}\otimes \Qp/\Zp)\ar[lld] \\
  & R^{1}\Crys_{1}(T_{p}A_{K}) \ar[r]
  & R^{1}\Crys_{1}(V_{p}A_{K}) \ar[r]
  &  \cdots
  \end{tikzcd}
\end{equation*}
attached to $T_{p}A_{K}$ using \Cref{les-exists}. The kernel of \[f: R^{1}\Crys_{1}(T_{p}A_{K})
  \longto R^{1}\Crys_{1}(V_{p}A_{K})\]
is exactly the right-hand side of \eqref{eq:Phi-ito-crys}, which is torsion, so we have
\[\Phi_{A_{K}}[p^{\oo}]\subset (R^{1}\Crys_{1}(T_{p}A_{K}))_{\tors}.\]
On the other hand, $R^{1}\Crys_{1}(V_{p}A_{K})$ is torsion-free (\Cref{crys-of-Qp-vs-tors-free}), so the kernel of $f$ contains $R^{1}\Crys_{1}(T_{p}A_{K})_{\tors}$. This completes the proof.
  \qed}
\end{theorem}


\bibliographystyle{abbrv}
\bibliography{biblio}

\end{document}